%% file: toms_2018.tex
\documentclass[manuscript, screen]{acmart}
\usepackage{booktabs} 

\usepackage{graphicx}
\usepackage{subfigure}
\usepackage{rotating}

\usepackage[ruled]{algorithm2e} 

\newcommand{\LL}{\mathcal L}

\newcommand\reals{{{\rm l} \kern -.15em {\rm R} }}
\newcommand{\ignore}[1]{}

\acmJournal{TOMS}
\acmVolume{9}
\acmNumber{4}
\acmArticle{39}
\acmYear{2017}
\acmMonth{3}

\setcopyright{acmcopyright}

\acmDOI{0000001.0000001}

\begin{document}
\title{Analysis and Performance Evaluation of Adjoint-Guided Adaptive Mesh
Refinement for Linear Hyperbolic PDEs Using Clawpack} 

\author{Brisa N Davis}
\affiliation{%
  \institution{University of Washington}
  \streetaddress{1410 NE Campus Parkway}
  \city{Seattle}
  \state{WA}
  \postcode{98195}
  \country{USA}}
\author{Randall J LeVeque}
\affiliation{%
  \institution{University of Washington}
  \streetaddress{1410 NE Campus Parkway}
  \city{Seattle}
  \state{WA}
  \postcode{98195}
  \country{USA}}

\renewcommand\shortauthors{Davis, B.N. et al.}
\renewcommand\shorttitle{Analysis and Performance Evaluation of
Adjoint-Guided Adaptive Mesh Refinement}

\begin{abstract}
Adaptive mesh refinement (AMR) is often used when solving time-dependent partial
differential equations using numerical methods. It enables time-varying regions of 
much higher resolution, which can be used to track discontinuities in the solution by selectively 
refining around those areas. The open source Clawpack software 
implements block-structured AMR to refine around propagating waves in the AMRClaw package. For problems 
where the solution must be computed over a large domain but is only of interest
in a small area this approach often refines waves that will 
not impact the target area. We seek a method that enables the 
identification and refinement of only the waves that will 
influence the target area.

Here we show that solving the time-dependent adjoint equation
and using a suitable inner product allows for a more precise 
refinement of the relevant waves. We present the adjoint methodology 
in general, and give details on how this 
method has been implemented in AMRClaw. Examples 
for linear acoustics equations are presented, 
and a computational performance analysis is conducted. 
The adjoint method is compared to AMR methods already available in the 
AMRClaw software, and the resulting advantages and disadvantages are discussed.
The code for the examples presented is archived on Github.
\end{abstract}

%
%


%
%


\thanks{This work is supported by the National Science Foundation,
  under an NSF Graduate Research Fellowship DGE-1256082 and grants
DMS-1216732 and EAR-1331412, as 
  well as the Graduate Presidential Fellowship, GO-MAP, University of Washington.}

\maketitle

\input{body-toms2018}

\end{document}

%% file: body-toms2018.tex
\section{Introduction}\label{sec:intro}
\indent Hyperbolic systems of partial differential equations
appear in the study of numerous physical phenomena involving wave
propagation. Methods for numerically calculating solutions to
these systems of PDEs have broad applications in
many disciplines. Complicating the development of numerical methods for
solving these systems is the fact that solutions often
contain discontinuities or localized steep gradients.
A variety of adaptive mesh refinement (AMR) techniques have been developed
to allow the use of much finer grids around discontinuities or localized
regions needing higher resolution, which generally propagate as the solution
evolves. 
Nonlinear hyperbolic problems that produce shock waves often require AMR 
in regions that can be easily identified by computing the local gradient. 
But AMR is often also extremely useful for linear hyperbolic equations with smooth 
solutions, particularly when the solution is of interest over
some small region relative to the size of the computational domain. 
Examples motivating this work include modeling tsunami propagation over the
ocean (for which linearized shallow water equations are suitable until the
waves reach shore) and earthquake modeling, where linear elasticity
equations are generally used. In both cases spatially varying
coefficients in the PDE lead to scattering and reflections that can
make it difficult to predict which portions of the domain must be
refined at any given time in order to capture the waves that will
ultimately reach the location of interest.  

We present a general approach to using the adjoint
equation to efficiently guide adaptive refinement in this situation, and
describe an implementation of this approach in the Clawpack software
\cite{CLAWPACK,mandli2016clawpack}.
This open source software has been developed since 1994 and can be applied to 
almost any hyperbolic PDE in 1, 2, or 3 space dimensions
by providing a Riemann solver, the basic building
block of the high-resolution Godunov-type finite volume methods that are
employed \cite{BaleLevMitRoss02,LangsethLeVeque00,LeVeque1997,Leveque1}.  
The AMRClaw package of Clawpack uses block structured AMR as developed in 
\cite{BergerOliger1984,BergerColella1989} and adapted to the
wave-propagation algorithms used in Clawpack in \cite{Berger1998}.
The AMR software is also used in the GeoClaw
variant of Clawpack developed
for modeling tsunamis, storm surge, and other geophysical flows (e.g.
\cite{BergerGeorgeLeVequeMandli:awr11,LeVequeGeorgeBerger2011,MandliDawson2014}).

For a hyperbolic PDE, the adjoint equation is another (closely related)
hyperbolic PDE that must be solved backwards in time,
starting with a disturbance at the location of interest. 
Waves propagating outward in this adjoint solution indicate the portion of
the domain that can affect the solution at the location and time of
interest. Taking a suitable inner product between the forward solution (or
an estimate of the local error in the forward solution) and the adjoint solution
allows us to determine what regions need to be refined at each earlier time.

We concentrate on linear (variable coefficient) problems in this paper
because in this case the adjoint equation is independent of the forward
solution, and can be computed {\em a priori} before solving the forward
problem.  The adjoint problem is generally solved on a much coarser grid
than will be used for the forward problem, so that this adds relatively
little computational expense.  For a nonlinear hyperbolic PDE an adjoint
problem can still be defined, but involves linearizing about a particular
solution of the forward problem. For nonlinear problems,
making full use of the adjoint to
guide adaptive refinement might require iterating between approximations to
the forward and adjoint problem as the forward problem is better resolved
in appropriate regions.  The work reported in this paper could form the
basis for such an extension, and has already proved very useful in its own
right for the applications mentioned above.
We also assume that the coefficients vary in space but not in time.
Time-varying coefficients could be handled in much the same way if the
focus is on a single output time of interest, but in most practical cases
the coefficients are time-invariant and hence the equation is autonomous in
time.  This makes it easy to extend the approach to handle problems where the
solution at one spatial location is of interest over a range of times, which
is the situation in the examples mentioned above, for example. 
The adjoint solution needs to be computed only once over a sufficiently long
time period, and snapshots in time of this one solution can be used to guide
the adaptive refinement regardless of the time period of interest.

In two space dimensions (with obvious modifications for 1 or 3 dimensions),
a linear hyperbolic system with coefficients varying
in space typically takes one of two forms, either the conservative form
\begin{equation}\label{eq:hypsys_conserv}
q_t(x,y,t) + (A(x,y)q(x,y,t))_x + (B(x,y)q(x,y,t))_y = 0,
\end{equation} 
or the non-conservative form
\begin{equation}\label{eq:hypsys_nonconserv}
q_t(x,y,t) + A(x,y)q_x(x,y,t) + B(x,y)q_y(x,y,t) = 0.
\end{equation} 
Of course \cref{eq:hypsys_nonconserv} differs from \cref{eq:hypsys_conserv} only if the coefficient
matrices vary in space.
Often the same physical system can be modeled by equations of either form,
depending on what variables are chosen to make up the vector $q(x,y,t)$.  For
example, linear acoustics can be written in the form
\cref{eq:hypsys_conserv} as a system of 3 equations where 
$q$ models small disturbances of density and
momenta in the $x$- and $y$-directions.  Alternatively the same acoustics
problem can be modeled 
in the form \cref{eq:hypsys_nonconserv} (with different coefficient
matrices) if instead the system
is written in terms of pressure and velocities.  Clawpack can be used to solve
equations in either form.  This is important to note since if the original
equation is in one form then the adjoint equation takes the other form.
As derived below in \cref{sec:adjeqn}, the adjoint equations arise from an integration
by parts and hence if the original equation is \cref{eq:hypsys_conserv} then
the adjoint equation is
\begin{equation}\label{eq:adjoint_nonconserv}
\hat q_t(x,y,t) + A(x,y)^T \hat q_x(x,y,t) + B(x,y)^T \hat q_y(x,y,t) = 0,
\end{equation} 
while if the original equation is \cref{eq:hypsys_nonconserv} then
the adjoint equation is
\begin{equation}\label{eq:adjoint_conserv}
\hat q_t(x,y,t) + (A(x,y)^T \hat q(x,y,t))_x + (B(x,y)^T \hat q(x,y,t))_y = 0.
\end{equation}

In this paper we first consider one-dimensional problems (one space
dimension plus time) for ease of exposition and because it is informative
to view the solutions in the $x$--$t$ plane in this case. The ideas extend
directly to more space dimensions and we also consider two-dimensional
acoustics in \cref{sec:2d_acoustics_example}. While not done in the work, 
the ideas presented here also extend directly to three dimensions.

A key component of any AMR algorithm is the criterion for deciding which grid cells
should be refined. 
Before this work, four different criteria were available in AMRClaw and GeoClaw for 
flagging cells that need refinement from one level to the next finer level. These are:
\begin{itemize}
\item a Richardson extrapolation
error estimation procedure
that compares the solution on the existing grid with the solution on a coarser
grid and refines cells where this error estimate is greater than a specified tolerance,
\item refining cells where the gradient of the solution (or an undivided difference of neighboring cell
values) is large,
\item refining cells where the surface elevation differs significantly from sea level 
(this refinement criteria is specific to GeoClaw for tsunami or storm surge modeling), or
\item some other user-specified criterion that examines the current solution locally.
\end{itemize}
In general these approaches will flag cells for refinement
anywhere that a specified refinement tolerance is exceeded,
irrespective of the fact that the area of interest may be only a subregion of
the full solution domain. 
To address this, recent versions of AMRClaw and GeoClaw also allow specifying
``refinement regions,'' space-time subsets of the computational domain where
refinement above a certain level can be either required or forbidden.  This
is essential in many GeoClaw applications where only a small region along the
coast (some community of interest) must be refined down to a very fine
resolution (often $1/3$ arcsecond, less than 10 meters) 
as part of an ocean-scale simulation. In this work we focus on 
acoustics examples, and in that context the region of interest might 
be dictated by the placement of a pressure 
gauge or other measurement device. 
These AMR regions can also be used to ensure that the simulation only refines  
around the waves of interest. However, since the regions are user-specified, placing 
them optimally often requires multiple attempts and careful examination of how 
the solution is behaving. This generally necessitates the use of coarser grid 
runs for guidance, which adds computational and user time requirements. 
This manual guiding of AMR may also fail to capture some waves that are 
important. For example, a user may determine, based on a coarser grid run,
that a portion of a wave is heading away from the region of interest 
and therefore forbid excessive refinement in that area in an effort 
to save computational time. However, this 
portion of the wave may later reflect off a distant 
boundary or material heterogeneity, 
causing it to have an unexpected impact on the region of interest. 

For any problem where only a particular area of the total solution is of
interest, a method that allows specifically targeting and refining the
grid in regions that influence this area of interest would allow for 
computational savings while also ensuring that the accuracy of 
the solution is preserved. 
This challenge is the motivation for our work presented here.


When solving the original time-dependent PDE,
which we will call the ``forward problem'', generally we have 
initial data specified at some initial time $t_0$
and the problem is solved forward in time to find the solution at the target
location at some later time $t_f$ (or, more typically, over some time range
of interest).
The adjoint equation (derived in \cref{sec:adjeqn}) must then be
solved backwards in time from the final time
$t_f$ to the initial time $t_0$.
The ``initial '' data for the adjoint equation, which must be specified at the final time, 
is non-zero at the target location and then spreads out into
waves as the adjoint equation is solved. The key idea is
that at any intermediate ``regridding'' time $t_r$, 
between $t_0$ and $t_f$, the only regions in space
where the forward solution could possibly affect the target location at time
$t_f$ are regions where the adjoint solution is nonzero. Moreover, by
computing a suitable inner product of the forward and adjoint solutions at
time $t_r$ it is possible to determine whether  the forward solution at
a given spatial point will actually impact the target location 
at the final time $t_f$, or whether it
can be safely ignored.  
Alternatively, by computing a suitable inner product 
of an error estimate in the forward solution and the adjoint solution at 
time $t_r$, it is possible to estimate how much the error at a given spatial 
point will affect the final accuracy at the target location.
This information can then be used to
decide whether or not to refine this spatial location in the forward solution.


\subsection{Previous Work in the Literature}

The idea of using adjoint equations to guide computations has been around 
for many years. For an overview of the various fields in which adjoints 
have been utilized see \cite{DavisLeVeque:adjoint2016}. 
Historically the adjoint method has been commonly used to 
improve the accuracy of a solution by calculating a correction 
term and adding it back into the solution, see for example the 
work of \cite{PierceGiles2000}. 
Adjoint equations have also been used to compute the sensitivity 
of the solution to changes in the input data. 
This has led to the adjoint equations being utilized for system control 
in a wide variety of applications such as shallow-water wave control 
\cite{SandersKatopodes2000},  
optimal control of free boundary problems \cite{Marburger2012}, and 
aerodynamics design optimization \cite{Jameson1988}. 
More relevant to this current work, the adjoint method has been used in
literature to guide adaptive mesh refinement.
\cite{BeckerRannacher2001} used  
the adjoint equation to relate the global error to errors in 
the physical quantities of interest for the particular problem 
being considered, and then used these \textit{a posteriori} error 
estimates to refine or coarsen the mesh. 
\cite{VendittiDarmofal2002} and \cite{VendittiDarmofal2003} 
have applied the same methods to 
compressible two-dimensional inviscid
and viscous flow problems. \cite{Park2004} extended the 
methods of \cite{VendittiDarmofal2002} to three-dimensional 
problems in the area of computational fluid dynamics. 
All of these, however, dealt with steady state problems. 

Some work on using the adjoint method to guide AMR for 
time-dependent problems has also been done within the 
finite volume community. This is a more expensive problem 
due to the fact that there is not a single time-invariant 
adjoint solution but rather the problem requires 
computing and saving adjoint snapshots at various times. 
Some examples include 
\cite{KastFidkowski2013} and \cite{LuoFidkowski2011} 
who use the adjoint method to guide AMR for both 
spatial and temporal grid refinement for the 
Navier-Stokes equations and \cite{Kouhi2015} 
who applies this method to the compressible 
Euler equations. All of these examples 
of using the adjoint method to guide AMR work with the 
discrete adjoint formulation, whereas in this work we 
derive the continuous adjoint and then discretize this hyperbolic 
system using the same finite volume methods as 
used for the forward problem.   
We presented the novel use of the continuous adjoint method in the 
context of the shallow water equations for tsunami modeling, as well as using 
the adjoint method to guide adaptive mesh refinement 
for problems in that field, in \cite{DavisLeVeque:adjoint2016}. 
Since then, \cite{Lacasta2018} 
has published work also using the continuous adjoint 
equations for the shallow water equations --- which they used 
to control their internal boundary conditions. 
For a comprehensive 
comparison of the continuous and discrete adjoint formulations
 see \cite{Jameson2000}.

\subsection{Outline}
We first introduce various methods for determining where to apply 
AMR in \cref{sec:flagging}, some of which have been implemented 
in AMRClaw and GeoClaw already and some of which were newly 
implemented in the Clawpack framework for this work.
We then present more detail on the adjoint equation and 
how it can be used to guide AMR in \cref{sec:adjeqn}, and some 
specifics on implementing this new method in the context 
of Clawpack are presented in \cref{sec:algorithm}. 
This basic method is expanded to take into account the 
approximated error in the calculated solution at each 
time step in \cref{sec:adjErrorFlag}. 
Examples 
of this method applied to the one- and two-dimensional 
variable coefficient linear acoustics equations 
are given in \cref{sec:1d_acoustics_examples,sec:2d_acoustics_example}, respectively, as well as an analysis 
of the computational performance of this method relative to 
the previously available AMR methods in Clawpack.
For completeness, the appendix contains a discussion of how to solve
the Riemann problem for the one-dimensional variable-coefficient
acoustics problem and its adjoint.

\section{Adaptive Mesh Refinement in Clawpack}\label{sec:flagging}
The block-structured mesh refinement used in AMRClaw consists of a set of 
logically rectangular grid patches at multiple levels of 
refinement. The coarsest level contains grids that cover the entire domain, 
with subsequent levels representing progressively finer mesh resolutions. 
Generally each level is refined in both time and space to preserve the
stability of the explicit finite volume method. Each level, other than the coarsest 
level, is properly nested within the grids that comprise the next coarsest level. 
For each time step on Level $L$ starting with $L=1$ for the coarsest level
(and recursively applying to the finest level
$L_{\max}$), the following steps are performed:
\begin{enumerate}
\item If $L>1$ then fill ghost cells around each patch at this level (from
neighboring Level $L$ patches, or by interpolating from Level $L-1$).
\item Take a time step of length $\Delta t_L$ appropriate to this level.
\item If $L<L_{\max}$, then
\begin{enumerate} 
\item If it is time to regrid, flag cells on Level $L$ that need refining,
and cluster these into rectangular patches to define the new Level $L+1$ grid
patches.
\item Take $\Delta t_L / \Delta t_{L+1}$ time steps on all grid patches at
Level $L+1$ by recursively applying this procedure.  
\item Update the Level $L$ solution in regions covered by Level $L+1$
patches by cell-averaging the finer grid solution, and in cells directly
neighboring  these patches as needed to preserve conservation.
\end{enumerate} 
\end{enumerate}

We assume the refinement ratio $\Delta t_L / \Delta t_{L+1}$ is an integer,
generally equal to the refinement ratio in space from Level $L$ to $L+1$ in
order to preserve stability, since Clawpack uses
explicit finite volume methods that require the Courant number be bounded by 1.  
See \cite{Berger1998} for more details on the steps above, and
\cite{BergerOliger1984,BergerColella1989} for general background on this
approach.  In this paper we
are concerned only with step 3(a),  and in particular the manner in which cells
at one level are flagged for refinement to the next level.
The clustering is then done using an algorithm of Berger and Rigoutsos
presented in \cite{BergerRigoutsos1991}, which attempts to limit the number of
unflagged cells contained in the rectangular patches
while also not introducing too many separate patches.
Some buffering is also done so that the refined patches extend a few grid cells
out from the flagged cells.  Waves can travel at most one grid
cell per time step, so this ensures that waves needing refinement will not
escape from the refined patches in a few time steps, but in general for a wave
propagation problem it is necessary to regrid every few time steps on each
level.  The regridding interval and width of the buffer zone are generally
related, and can each be set as input parameters in the code.
Because of the need for frequent regridding, and the desire to minimize the
number of needlessly refined cells, the methodology utilized for
selecting which cells to flag will have a significant impact on both the
accuracy of the results and the time required for the computation to run.

\subsection{Flagging Methods Available in Clawpack}\label{sec:flaggingmeths}
Prior to this work there were
two built-in flagging methods available in AMRClaw, undivided differences and 
Richardson extrapolation. Both of these methods flag a cell if the 
quantity being evaluated is greater than some specified tolerance. 
A third flagging method, specific to GeoClaw, determines whether to 
flag a cell based on the surface elevation of the water.
For all of the methods of refinement, if the user has specified 
limitations on certain regions of the domain (either requiring 
or forbidding flagging to occur) then these limitations are also 
taken into account at the flagging step in the code.

{\em Undivided Differences. }
This is the default flagging routine that is used in AMRClaw. This 
routine evaluates the maximum max-norm of the undivided differences
between a given grid cell and its two, four, or six neighbors in 
one, two, or three space dimensions respectively.
 If this maximum 
max-norm is greater than the specified tolerance then the cell 
is flagged for refinement. Therefore, this flagging method refines 
the mesh wherever the solution is not sufficiently smooth. 
This approach often works very well for problems with shock 
waves where the goal is to refine around all shocks. 
In this work, we will refer to this flagging method as 
\textit{difference-flagging}.

{\em Richardson Extrapolation. }
This second approach to flagging is based on using Richardson extrapolation 
to estimate the error in each cell. For each grid on Level $L$,
this is done by 
\begin{itemize}
\item advancing the current grid a time step, as normal, 
\item taking an extra time step on the current grid,
\item generating a coarsened grid (by a factor of 2 in each direction)
and taking one time step on this grid,
\item comparing the solution on these two grids to estimate the 
one-step error introduced in a single time step with the current Level $L$
mesh size.  Cells are flagged for refinement to Level $L+1$ where this
estimate is above a given tolerance.
\end{itemize}
Note that this approach requires one additional time step on the Level $L$
grid and one time step on the coarsened grid, relative to difference-flagging.
Therefore, it is more expensive than difference-flagging, but
has the advantage of refining based on estimates of the error in the solution 
rather than simply anywhere that the solution is not sufficiently smooth.
In this work we will refer to this flagging method as \textit{error-flagging}.
 

\ignore{
\subsubsection{Surface Elevation}
This flagging method is specific to GeoClaw, and is in fact the 
default for that code. It evaluates the surface height of the water 
in each grid cell (the sum of water depth and underlying topography), 
and flags the cell if the difference between the 
surface height and sea level for the ocean at rest is above some tolerance. 
In this work we will refer to this flagging method as \textit{surface-flagging}.
As described in the previous section, 
cells in GeoClaw can also be flagged if they are in a region that 
enforces refinement. Defining regions where refinement is either required 
or not allowed is common for GeoClaw simulations, due to the large 
total simulation area (e.g., the Pacific Ocean)
and the level of resolution required to 
accurately model small regions of interest (e.g., one harbor at 10 m
resolution). In addition, another option in
GeoClaw allows flagging where the fluid speed is above 
a given speed tolerance.
}

\subsection{Flagging Using the Adjoint Method}\label{sec:adjoint}
The adjoint-flagging method is the focus of this work. Here a very brief 
description of how the method can be used for AMR is given, 
to provide a general framework for understanding our goal for the 
adjoint equations. 
 \Cref{sec:adjeqn} outlines the derivation of the equations 
required for this method, and \cref{sec:algorithm} 
contains a detailed description of the algorithms used to implement this 
method. \Cref{sec:adjErrorFlag} expands the basic adjoint method to take 
into account the approximated error in the calculated solution 
of the forward problem at each time step. 
The end goal of this work is to determine which regions 
of the domain need to be refined at any given time
by using the adjoint solution to determine what  
portions of the current forward solution 
will affect the target location.

The steps needed to use this 
flagging method are the following:
\begin{enumerate}
\item Determine the appropriate adjoint problem based on the forward 
problem.
\item Solve the adjoint problem backward in time, saving snapshots of the solution at 
various times.
\item When solving the forward problem, at each regridding time:
\begin{enumerate}
\item identify the snapshot(s) of the adjoint solution that bracket the
current time in the forward solution,
\item find the value of each of the identified adjoint snapshots at the spatial location of  
each grid cell by interpolation,
\item take a suitable inner product between the forward solution and each of the 
identified adjoint snapshots, 
\item and flag each cell if the maximum of these inner products 
for that grid cell is above a certain tolerance.
\end{enumerate}
\end{enumerate}
Note that there are several extra steps required for this flagging method, 
which we will refer to as \textit{adjoint-flagging}, when 
compared to the flagging methods currently available in AMRClaw. 
However, much of this has been automated in the work described here, and if the
adjoint problem is solved on a relatively coarse grid then the computational
time added by the adjoint solution and inner products may be small relative
to the time potentially saved by refining in a more optimal manner. 
Also note that, rather than considering the inner product between 
the forward problem and various selected adjoint snapshots, we could 
instead start from the appropriate adjoint snapshot and take a few small 
time steps to get to the time corresponding to the current regridding time 
for the forward problem. Since we are trying to minimize the amount 
of work required while solving the forward problem, we have selected to go 
with the method presented above --- as it requires less computations. 
However, even though the method we have selected is less computationally 
intensive, it has proven to be extremely effective in enabling us to 
guide AMR for the forward problem. 

Finding the appropriate adjoint problem analytically and 
implementing the adjoint solver is a one-time cost 
for each type for forward problem. For example, once the 
appropriate adjoint problem is found for the acoustics equations, any 
acoustics forward problem can use the same adjoint problem. 
Note that the adjoint equation is also a linear hyperbolic equation, and can therefore 
be solved using the same software as the forward problem. 
To see how this can be accomplished see  
\cref{sec:1d_acoustics_examples,sec:2d_acoustics_example} 
where the adjoint problems are found for one- 
and two-dimensional linear acoustics 
problems. 

When using adjoint-flagging we consider 
two different options when it comes to taking the inner product, 
stemming from a slightly different refinement criterion analogous to the
difference-flagging vs.\ error-flagging described above. 

{\em Considering the forward solution.}
The first option focuses on the magnitude of the forward solution, and asks 
the question:
\begin{description}
\item At a given regridding time $t_r$, what portions of the forward solution 
will eventually affect the target location?
\end{description}
To answer this question we take the inner product between the 
forward solution and the adjoint solution at the appropriate complementary  
time and flag the cells where this inner product is above 
the given tolerance. 
This has the advantage of flagging only cells that contain 
information that is pertinent to the target location and time 
interval, which is 
not possible with any of the flagging methods currently available 
in AMRClaw and GeoClaw.
In this work we will refer to this flagging method as 
\textit{adjoint-magnitude flagging}, which we describe in detail in
\cref{sec:adjeqn,sec:algorithm}.
While this is often quite effective, it can be difficult to choose a suitable
magnitude tolerance for flagging cells (relative to the desired accuracy of the
solution).  Moreover, regions where the inner product is largest do not
necessarily contribute the most error to the final computed solution at the
location of interest.
Whether these regions actually require more refinement depends also on
the smoothness of the solution and hence the accuracy of the finite volume 
solution.

{\em Considering the error in the forward solution.}
The second option focuses on estimating the error in the forward solution, 
and asks the question:
\begin{description}
\item At a given regridding time $t_r$, what portions of the forward solution 
will introduce a significant amount of error that will eventually affect 
the target location?
\end{description}
This is the question we really want to answer in deciding where to refine,
but requires more work.
To answer this question we estimate the one-step error in the forward solution 
using the Richardson extrapolation algorithm of AMRClaw, take the inner product 
between this error estimation and the adjoint solution at the 
appropriate complementary time, and flag the cells where this inner product 
is above a given tolerance. 
This method still has the advantage of flagging only cells that contain 
information that will reach the target location, but also only if the error
estimate indicates that the error in this part of the solution is signficant
on the current grid resolution.
However, this flagging method is slower than 
adjoint-magnitude flagging because it requires the 
estimation of the one-step error in each grid cell using Richardson 
extrapolation. We will refer to this flagging 
method as \textit{adjoint-error flagging}, and consider it further in
\cref{sec:adjErrorFlag} after describing the basic ideas and algorithms in
the simpler context of adjoint-magnitude flagging.

\section{The Adjoint Equation and Adjoint-Magnitude Flagging}\label{sec:adjeqn}
For readers not familiar with the concept of an adjoint equation, 
\cite{DavisLeVeque:adjoint2016} contains a description of the 
basic idea in the context of an algebraic system of 
equations that may be beneficial 
in appreciating this method, along with details of the adjoint AMR procedure
and some motivating examples in the context of tsunami modeling.  
Here, we will present the mathematics 
behind finding the adjoint equation in the general context of time-dependent 
linear hyperbolic equations, and 
\cref{sec:algorithm} describes the algorithms necessary to 
implement this method in AMRClaw and GeoClaw.

Suppose $q(x,t)$ is the solution to the time-dependent linear equation 
(with spatially varying coefficients)
\begin{equation} \label{eq:qtAqx}
q_{t}(x,t) + A(x)q_{x}(x,t) = 0, \quad a\leq x \leq b, \quad t_0\leq t \leq
t_f
\end{equation} 
subject to some known initial conditions, 
$q(x,t_0)$, and some boundary conditions at $x=a$ and $x=b$. 
Here $q(x,t) \in \reals^m$ for a system of $m$ equations and
we assume $A(x) \in \reals^{m\times m}$ 
is diagonalizable with real eigenvalues at each $x$, so that
\cref{eq:qtAqx} is a hyperbolic system of equations.
For the discussion and examples presented in this work 
we assume that we start with an equation in the form
\cref{eq:hypsys_nonconserv}, but we could equally well start with an equation
in the form \cref{eq:hypsys_conserv}.

Now suppose we do not care about the solution everywhere, but only about
the value of a linear functional
\begin{equation}\label{eq:J_general}
J = \int_a^b \varphi^T (x) q(x,t_f) dx
\end{equation} 
for some given $\varphi (x)$ at the final time (or, more typically, over a
range of times as considered below).
Although any function $\varphi (x)$ could be specified in the method and
software developed here, the case we will focus on in the examples is the
common situation where we are only interested in the solution in one small
spatial region (e.g., one coastal community in the case of tsunami
modeling), but we need to solve the equations over a much larger region
(e.g., the entire ocean) in order to determine the solution in this region.
In this case it is natural to define $J$ by choosing $\varphi (x)$ to be a
delta function centered at the \textit{point of interest}, 
$\varphi (x) = \delta (x - x_p)$.  
Or more realistically (and better computationally), we can take $\varphi
(x)$ to be a box function or Gaussian centered about $x_p$ with mass 1,
corresponding to $J$ being a spatial average of $q(x,t_f)$ near the location
of interest.  

If $\hat{q}(x,t)\in\reals^m$ is any other function then taking the inner
product of $\hat q$ with 
\cref{eq:qtAqx} and integrating over the space-time domain yields
\begin{equation}\label{eq:q_integral}
\int_{t_0}^{t_f}\int_a^b 
\hat{q}^T(x,t)\left(q_{t}(x,t)+A(x)q_{x}(x,t)\right)dx\,dt = 0.
\end{equation}
Then integrating by parts first in space and then in time yields the equation 
\begin{equation}\label{eq:intbyparts}
\int_a^b   \left.\hat{q}^Tq\right|^{t_f}_{t_0}dx 
+ \int_{t_0}^{t_f} \left.\hat{q}^TAq\right|^{b}_{a}dt 
- \int_{t_0}^{t_f} \int_a^b  q^T\left(\hat{q}_{t} +
\left(A^T\hat{q}\right)_{x}\right)dx\,dt = 0.
\end{equation} 
By choosing $\hat q(x,t)$ to satisfiy the adjoint equation, 
\begin{equation}\label{adjoint1}
\hat{q}_{t}(x,t) + (A^T(x)\hat{q}(x,t))_{x} = 0,
\end{equation}
solved backward in time from  $\hat{q}(x,t_f) = \varphi (x)$,
and selecting the appropriate boundary conditions for $\hat{q}(x,t)$ 
such that the second integral in \cref{eq:intbyparts}
vanishes (which varies based on the 
specific system being considered), we can eliminate all
terms from \cref{eq:intbyparts} except the first term, to obtain 
\begin{equation}\label{eq:q_equality}
\int_a^b \hat{q}^T(x,t_f)q(x,t_f) dx = \int_a^b
\hat{q}^T(x,t_0)q(x,t_0)dx.
\end{equation} 
Therefore, the integral of the inner product
between $\hat{q}$ and $q$ at the final time is equal 
to the integral at the initial time $t_0$:
\begin{equation}
J = \int_a^b\hat{q}^T(x,t_0)q(x,t_0)dx.\label{eq:J}
\end{equation}
Note that we can replace $t_0$ in \cref{eq:q_integral} with any 
$t_r$ at which we wish to do regridding ($t_0 \leq t_r \leq t_f$), 
which would yield \cref{eq:J} with $t_0$ replaced by $t_r$. 
From this we observe that the locations where the inner product
$\hat{q}(x,t_r)^Tq(x,t_r)$ is large at the regridding time
are the areas that will have a
significant effect on the functional $J$.
These are the areas where the solution should be refined at time $t_r$ 
if we are using adjoint-magnitude flagging.
To make use of this, we must first solve the adjoint equation
\cref{adjoint1} for $\hat q(x,t)$.
This requires using 
``initial'' data $\hat{q}(x,t_f)=\varphi(x)$, so
the adjoint problem must be solved backward in
time. The strategy used for this is discussed in  
the next section.  

\section{Implementing the Adjoint Method}\label{sec:algorithm}

In principle the adjoint-flagging methodology could be used with 
any software that uses time-dependent AMR. In this work we are examining the 
specifics of implementing it in the context of AMRClaw. 
This section discusses the special considerations required for coding 
the adjoint method in Clawpack, and the resulting algorithms. 

\subsection{Solving the Adjoint Equation}

Consider the one dimensional problem \cref{eq:qtAqx} and recall that
the adjoint equation \cref{adjoint1} has the form
\begin{align*}
\hat{q}_t + \left( A^T(x) \hat{q}\right)_x  = 0,
\end{align*}
where the initial condition for $\hat{q}$ is given at the final time,
$\hat{q}(x,t_f) = \varphi (x)$, and is selected to highlight the
impact of the forward solution on some region of interest. 

Clawpack is designed to
solve equations forward in time, so slight modifications must be made to the 
adjoint problem to make it compatible with our software. 
Define the function
\begin{align*}
\tilde{q}(x,\tilde{t}) &\equiv \hat{q}(x,t_f - \tilde{t})
\end{align*}
for $\tilde{t} > 0$.
This function satisfies the \textit{modified adjoint equation} 
\begin{align}\label{eq:1d_adjoint_modified}
&\tilde{q}_{\tilde{t}} - \left( A^T(x) \tilde{q}\right)_x = 0 
&&x \in [a,b], \hspace{0.1in}\tilde{t} > 0\\
&\tilde{q}(a,\tilde{t}) = \hat{q}(a,t_f - \tilde{t})  &&0 \leq \tilde{t} \leq t_f - t_0 \nonumber\\
&\tilde{q}(b,\tilde{t}) = \hat{q}(b,t_f - \tilde{t}) &&0 \leq \tilde{t} \leq t_f - t_0 \nonumber
\end{align}
with initial condition $\tilde{q}(x,0) = \varphi (x)$ given at the initial
time $\tilde t_0 = 0$. This
problem is then solved using the Clawpack software. Snapshots of this solution
are saved at regular time intervals, $\tilde{t}_0 = 0,~\tilde{t}_1,~\cdots,~\tilde{t}_f = t_f - t_0$. 
Note that $\tilde{t}$ has the interpretation of time remaining to $t_f$, which 
will be useful when selecting which adjoint snapshots to consider at a given 
regridding time $t_r$. This is discussed in the next section. 

When solving the forward problem, the saved 
snapshots of the adjoint solution are retrieved from the appropriate 
output folder and the adjoint solution $\tilde q$ at each snapshot is 
saved in an `adjoints' data structure. This structure is then utilized 
throughout the code whenever we need to take the inner product between 
the adjoint solution and either the forward solution or the Richardson error 
estimate of the forward solution. 

The basic mechanism through which Clawpack solves hyperbolic problems is 
by using a Riemann solver to calculate the waves generated between 
each set of adjacent cells at any point in time. See \cite{Leveque1} for 
more details of the wave propagating algorithms that are used. Therefore, an 
appropriate Riemann solver is needed for any problem being solved using 
this software package. 
For linear systems the Riemann solution is generally easy to compute in terms
of the eigenstructure of the coefficient matrix, for either nonconservative
or conservative linear systems.  We provide details for both the acoustics
equation and its adjoint in the appendix, and the software implementation of
these solvers are available at \cite{riemannCode}.

\ignore{
Riemann solvers for a variety of problems 
are available with the basic Clawpack software \cite{CLAWPACK}, including the
acoustics equations in the form used here.
For the adjoint problems a new Riemann solvers was required which 
can be found at \cite{riemannCode}.

We have previously noted that since our forward problem 
is in nonconservation form the adjoint problem we achieve 
is in conservation form. 
It is however possible to begin with a forward problem in conservative 
form, which would yield an adjoint equation in nonconservation form. This 
is significant in the context of Clawpack due to the solution methods 
utilized for solving these two different forms of hyperbolic systems 
of equations. 
Let $Q(x_i,t_n)$ be the 
calculated solution that approximates $q(x_i,t_n)$ where $x_i$ is 
a location in the discretized spatial grid. Then for a problem in 
nonconservation form the Riemann problem can be solved by resolving 
the jump from $Q(x_{i-1},t_n)$ to $Q(x_i,t_n)$ into waves 
propagating into neighboring cells. 
When we have a problem in conservative form we can 
instead use the f-wave formulation as described in 
\cite{BaleLevMitRoss02} and
\cite{LeVeque2011}, which splits the flux difference 
$f\left(Q(x_{i-1},t_n)\right) - f\left(Q(x_i,t_n)\right)$ into 
propagating waves. In this work we are considering forward 
problems in nonconservation form, so 
the Riemann solvers we developed for the adjoint 
problems use the f-wave formulation.
}

\subsection{Using Adjoint Snapshots}

With the adjoint solution in hand, we now turn to solving the forward problem. 
During this solution process it is unlikely that solution data for the
adjoint problem will be available at all the times needed for regridding, nor
will it generally be available on as fine a grid as the forward solution,
since the adjoint solution is generally solved only on a coarse grid.
As refinement occurs in space for the forward problem, maintaining the stability
of the finite volume method requires that refinement must also occur in time,
which will further exacerbate the issue of adjoint solution data not 
being available at the needed spatial and temporal locations. 
To
address this issue, the solution for the adjoint problem at the necessary 
locations is approximated using  linear or bilinear interpolation from the data present
on the coarser grid for one- or two-dimensional problems, respectively.

Typically we are interested in a time range, between $t_s$ and $t_f$, 
at the location of interest rather than a single point in time. For example, 
if we were modeling an acoustic wave we could be interested in the time range 
between when the first and last waves reach our pressure gauge, or in a
tsunami simulation we are interested in accurately capturing the waves
reaching a tide gauge over some time range. See  \cite{DavisLeVeque:adjoint2016}
for more details and examples for this latter application.
This time range comes into play when we are 
considering which adjoint snapshots to take into account at 
each time step of the forward problem. Suppose that we are solving the 
forward problem from $t_0$ to $t_f$ and that we are 
currently at regridding time $t_r$ in the solution process. 
Since the time for our modified adjoint problem,
 $\tilde{t}$, has the interpretation of being 
the time remaining to $t_f$, when evaluating which portions 
of the current solution will 
affect the target location at $t_f$ we need to consider 
the modified adjoint at $\tilde{t} = t_f - t_r$. 
Similarly, when 
evaluating which portions of the current solution will 
affect the target location at $t_s$ we need to consider 
the adjoint at $\tilde{t} = t_s - t_r$ (using the fact that the equations are
autonomous in time). Because we are 
actually interested in evaluating which portions of 
the current solution will affect the target location 
in the time range between $t_s$ and $t_f$, we need 
to consider the adjoint solutions snapshots for which  
$t_s - t_r \leq \tilde{t} \leq t_f - t_r$. This results 
in a list of adjoint snapshots, 
$\tilde{t}_m$ for $m = m_1, m_2, \cdots, M$ that need to be considered 
for each time $t_r$ that we wish to perform flagging for 
the forward problem. 
To be conservative, we expand the list to include $\tilde{t}_{m_1 - 1}$ 
(if $\tilde{t}_{m_1} > 0$) and $\tilde{t}_{M + 1}$ (if $\tilde{t}_M < t_f - t_0$). 

We then take the inner product between the forward solution 
and each adjoint snapshot on our list. 
(In \cref{sec:adjErrorFlag} we modify this to use the one-step
error estimate rather than the forward solution itself.)
Since we are 
concerned with flagging the cell if any of these inner products exceeds 
the specified tolerance, we just keep track of the maximum inner product calculated. 
If this value exceeds the tolerance, the cell is flagged. 
Algorithm \ref{alg:flagging} describes in pseudo-code the basic 
flagging procedure, although in the actual code base this functionality 
is spread throughout various files. 

\begin{algorithm}[t]
\SetAlgoNoLine
\KwIn{Forward solution on some level $L<L_{\max}$ at a particular regridding time $t_r$, and 
adjoint solution snapshots}
\KwOut{Flagged cells needing refinement to level $L+1$}
    $list$ = adjoint snapshots at times $\tilde{t}_{m_1}, \tilde{t}_{m_2}, \cdots\tilde{t}_{M}$ such that 
    $t_s - t_r \leq \tilde{t}_m \leq t_f - t_r$\\
    if $\tilde{t}_{m_1} > 0$ add adjoint snapshot at time $t_{m_{1}-1}$ to $list$\\
    if $\tilde{t}_M < t_f - t_0$ add adjoint snapshot at time $t_{M+1}$ to $list$\\
    \For{each grid patch \textnormal{g} at level \textnormal{L}}
        {
        initialize max inner product to 0 at each grid point\\
         \For{each adjoint snapshot in \textnormal{list}}
            {
            \For{each grid patch in adjoint snapshot}
                {
                 \If{the adjoint and forward patches overlap}
                 {
                 \For{each cell in grid patch \textnormal{g}}
                    {
                 interpolate the adjoint snapshot in space and time to cell center\\
                 if using adjoint-error flagging: estimate error in forward solution (\cref{sec:adjErrorFlag}) \\
                 calculate appropriate inner product and save if greater than current recorded maximum
                 }
                }
              }
            }
        \For{each cell on grid patch \textnormal{g}}
            {
            if max inner product exceeds tolerance: flag cell
           }
        }
\caption{Flagging Cells For Refinement}
\label{alg:flagging}
\end{algorithm}

As in the standard
AMRClaw and GeoClaw codes, the user must choose a tolerance and some
experimentation may be required to choose a suitable tolerance, related to
the scaling of the solution.   This form of
adjoint-magnitude flagging allows the code to avoid refining regions of the
solution that cannot possibly affect $J$ 
(the inner product will be indentically zero in such regions). 
Moreover, the inner product $|\hat q^T(x,t_r)q(x,t_r)|$
can be viewed as a measure of the sensitivity of
$J$ to changes to the solution near $x$ at the regridding time, and hence
regions where this is large may be important to refine.  However, just
because this is large at some $x$ does not necessarily mean that the solution
is inaccurate at this location and needs refinement.
The next section explores the extension of this approach to incorporate error
estimation.

\section{Adjoint-Error Flagging}\label{sec:adjErrorFlag}
To extend the adjoint approach to use adjoint-error flagging,
we need to consider the errors being made at each time 
step. Let $Q(x,t_n)$ 
be the calculated solution at $t_n$ that approximates $q(x,t_n)$. 
Then the error in the functional $J$ at the final time is given by 
\begin{equation*}
\int \hat{q}^T(x,t_f) \left[Q(x,t_f) - q(x,t_f)\right] dx.
\end{equation*}
Recall that the PDE we solve is assumed to be linear and autonomous in time,
and let $\LL(\Delta t)$ represent the
solution operator over time $\Delta t$, so we can write
$q(\cdot, t+\Delta t) = \LL(\Delta t)q(\cdot, t)$.  
For notational convenience we will also write this as 
$q(x, t+\Delta t) = \LL(\Delta t)q(x, t)$,  
but note that $\LL(\Delta t)$ cannot be applied pointwise, so this is really
shorthand for 
$q(x, t+\Delta t) = [\LL(\Delta t)q(\cdot, t)](x)$.  
Also note that $\LL$ satisfies the semigroup property, $\LL(t_2-t_1)\LL(t_1-t_0) =
\LL(t_2-t_0)$.

Now let $\tau(x,t_n)$ be the one-step error at time $t_n$, defined as the error
that would be incurred by taking a single time step of length $\Delta t_n =
t_n-t_{n-1}$, starting with data
$Q(x, t_{n-1})$:
\begin{equation*}
\tau(x,t_n) = Q(x, t_n) - \LL(\Delta t_n)Q(x, t_{n-1}),
\end{equation*}
which is approximately $\Delta t_n$ times the local truncation error (LTE). 
Then the global error grows according to the recurrence 
\begin{equation*}
Q(x,t_n) - q(x,t_n) = \LL(\Delta t_n)(Q(x, t_{n-1})-q(x, t_{n-1})) + \tau_n,
\end{equation*}
which, together with the assumption of no error in the initial data 
yields, after $N$ steps,
\begin{equation}\label{eq:finalerror}
Q(x,t_N) - q(x,t_N) = \sum_{n = 1}^{N}\LL(t_N - t_n) \tau (x,t_n).
\end{equation}

Now note that \cref{eq:q_equality} implies that
$\int \hat{q}^T(x,t_N)\LL(t_N-t_n)q(x,t_n) dx = \int
\hat{q}^T(x,t_n)q(x,t_n)dx$.  This same expression applies if we replace
$q(x,t_n)$ by any other function of $x$, since the solution operator $\LL$ then
evolves this data according to the original PDE.  Hence, in particlar,
\begin{equation}\label{eq:erroreq}
\int \hat{q}^T(x,t_N)\LL(t_N - t_n)\tau(x,t_n) dx = \int \hat{q}^T(x,t_n)\tau (x,t_n) dx.
\end{equation}

Bearing this in mind, we turn back to considering the error in our 
functional of interest, given by 
\cref{eq:J_general}. Let $E_J$ be the error in our functional at the final 
time $t_N \equiv t_f$: 
\begin{equation*}
E_J = \int \hat{q}^T(x,t_N)Q(x,t_N) dx - \int \hat{q}^T(x,t_N)q(x,t_N)dx.
\end{equation*}
Using \cref{eq:finalerror,eq:erroreq} gives us 
\begin{equation}\label{eq:EJsum}
\begin{split} 
E_J &= \int \hat{q}^T(x,t_N)\left[Q(x,t_N) - q(x,t_N)\right]dx \\
&= \int \hat{q}^T(x,t_N)\sum_{n = 1}^{N}\LL(t_N - t_n) \tau (x,t_n)dx \\
&= \sum_{n = 1}^{N}\int \hat{q}^T(x,t_N)\LL(t_N - t_n) \tau (x,t_n)dx \\
&= \sum_{n = 1}^{N} \int \hat{q}^T(x,t_n)\tau (x,t_n)dx.
\end{split} 
\end{equation}
Therefore the error in our functional at the final time 
is equal to the sum of integrals of the inner product of the
adjoint and the one-step error at each time step. From this we can 
observe that the locations where the inner product $\hat{q}^T(x,t_n)\tau (x,t_n)$ 
is large at the time $t_n$ are the areas that will have a significant effect on 
the final error. 
Moreover we can attempt to use estimates of the one-step error to flag only
the cells that contribute excessively to the estimated error in $J$.

Now suppose we wish to limit the error $E_J$ to a maximum value of $\epsilon$.  
Then we want 
\begin{equation*}
\left| \sum_{n = 1}^{N} \int \hat{q}^T(x,t_n)\tau (x,t_n)dx\right| \leq \epsilon.
\end{equation*}
Note that
\begin{align}\label{eq:sumabs}
\left| \sum_{n = 1}^{N} \int \hat{q}^T(x,t_n)\tau (x,t_n)dx\right| 
&\leq \sum_{n = 1}^{N}\Delta t_n \left| \frac{1}{\Delta t_n}
  \int \hat{q}^T(x,t_n)\tau (x,t_n)dx\right|  \\
&\leq (t_N-t_0) \max\limits_{n}\frac{1}{\Delta t_n}  
\int \left| \hat{q}^T(x,t_n)\tau (x,t_n)\right|dx .
\end{align}
We can enforce this bound by requiring that
\begin{equation}\label{eq:tserror}
 \int \left|  \hat{q}^T(x,t_n)\tau (x,t_n)\right| dx \leq \epsilon
\Delta t_n / (t_N - t_0) \equiv \epsilon_n
\end{equation}
for each $n$.
In the AMRClaw software the  
Richardson error estimation procedure described in \cref{sec:flaggingmeths}
provides an estimate of the one-step error, and this can be used to
approximate the contribution to the error $E_J$ in the functional that
results from the step at $t_n$. 

Let $Q(x_i, t_n; L)$ represent the numerical solution in a grid cell at level
$L$, where $i$ ranges over an index set $i\in I_L$ giving the indices of cells
at this level that are not covered by finer grid patches.  
Only these cells are part of the best approximation to $q(x,t_n)$ that is
composed of the finest-level approximation available at each point in the domain.
Similarly, let
$\tau(x_i, t_n; L)$ be the estimated one-step error in such a cell, and let $\hat
Q(x_i, t_n; L)$ be the numerical approximation of the adjoint solution 
interpolated to $(x_i, t_n)$.  Then the error estimate \cref{eq:EJsum} for $E_J$
can be approximated by
\begin{equation}\label{eq:EJest}
E_J \approx E_J^n =  \sum_{L=1}^{L_{max}} \sum_{i\in I_L} \hat Q^T(x_i,t_n;L)
\tau(x_i,t_n;L) \, \Delta x_L.
\end{equation} 
In practice we choose the maximum refinement level $L_{max}$ in advance along
with some desired tolerance $\epsilon$ and hope to achieve an error in $J$
that is no larger than $\epsilon$.  This makes the implicit assumption that if
we refined everywhere to the finest level then we would be able to obtain at
least this much accuracy in $J$.  Our goal is to acheive $|E_J| < \epsilon$
without refining everywhere, by refining only the areas on each level $L$ where
the error estimate is too large.  Choosing the optimal refinement pattern would
be an unsolvable global optimization problem over all grids, so we must use some
heuristic to decide what is ``too large'' as we sequentially refine each level
to the next resolution.  We take the approach of trying to enforce
\begin{equation}\label{eq:EJlevel}
\left| \sum_{i\in I_L} \hat Q^T(x_i,t_n;L) \tau(x_i,t_n;L) \, \Delta x_L\right|
< \epsilon_n/L_{max}
\end{equation} 
on each level, so that the allowable error is equally distributed among levels.
This suggests a way to flag cells at level $L$ for refinement to level $L+1$: 
flag cells starting with the one
having the largest error estimate and continue flagging (i.e., removing cells
from the index set $I_L$)  until the sum of the
remaining estimates satisfies \cref{eq:EJlevel}.

This is still impossible to implement in the context of AMRClaw, where there may
be many patches at each level that are handled sequentially (or in parallel with
OpenMP), and instead we require a pointwise tolerance
$\epsilon_n^L$ so that we can simply flag any cell where 
\begin{equation}\label{eq:EJptwise}
\left|\hat Q^T(x_i,t_n;L) \tau(x_i,t_n;L) \, \Delta x_L\right| > \epsilon_n^L.
\end{equation} 
One simple choice is $\epsilon_n^L = (\epsilon_n/L_{max})/(b-a)$, 
where $(b-a)$ is the length of the domain,
since $\sum_{I_L} \Delta_x \leq (b-a)$. This
ignores the fact that in general only part of the domain is refined to level
$L$, and so this can be improved by
replacing $(b-a)$ by the total length of level $L$ grid patches.
This is the approach we use in the first time step $n=0$ to choose the initial
grid patches at each level.   
This still ignores the fact that in some cells the
error will be much smaller than the tolerance,
and hence in other cells a larger error should be allowed.
So at later regridding times we 
choose $\epsilon_n^L$ based on the error estimates from the
previous regridding step, assuming a similar distribution of errors (generally
valid since we regrid every few time steps on each level). 
We save all the error estimates as each grid patch is
processed.  At the next regridding time for each refinement level, $L$, 
we sort the error estimates 
from all of the grids corresponding to level $L$ from smallest to largest and sum them up until the
cumulative sum of the first $j$ sorted errors reaches $\epsilon_n/L_{max}$, and
then set $\epsilon_n^L$ to be the last included term,
\begin{equation}\label{eq:epsnL}
\epsilon_n^L = \left|\hat Q^T(x_j,t_n;L) \tau(x_j,t_n;L) \, \Delta x_L\right|
\end{equation} 
In the next regridding step we flag any cell on level $L$ for which the pointwise error
estimate exceeds this value.

\section{One-Dimensional Variable Coefficient Acoustics}
\label{sec:1d_acoustics_examples}
In the next two sections we give some examples of the performance of the
algorithms developed here relative to the older AMRClaw algorithms.
We begin with two one-dimensional linear 
acoustics examples for ease of visualization.
Our first example showcases the 
capabilities of adjoint-flagging in a rather simple setting, and the 
second showcases a more complex scenario. These two examples 
highlight the strengths of adjoint-flagging 
when compared to the flagging methods currently available in AMRClaw. 
For examples using adjoint-flagging in a two dimensional acoustics 
context see \cref{sec:2d_acoustics_example}. 
For some examples using adjoint-flagging in a two dimensional shallow 
water equations context see \cite{Borrero2015}, 
\cite{DavisLeVeque:adjoint2016}, and \cite{adjointCode}.

Consider the linear acoustics equations in one dimension in a piecewise
constant medium, 
\begin{equation}\label{eq:1d_linearacoustics}
\begin{split} 
p_t(x,t) + K(x) u_x(x,t) &= 0,\\
{\rho (x)}u_t(x,t) + {p_x(x,t)} &= 0,
\end{split}
\end{equation} 
in the domain $x \in [a,b], t > t_0$, with solid wall (reflecting)
boundary conditions at each boundary,
\begin{equation}\label{acou1bc}
u(a,t) = 0, \hspace{0.1in} u(b,t) = 0, \qquad t  \geq t_0. 
\end{equation}
Setting 
\begin{align}\label{eq:Aacoustics}
A(x) = \left[ \begin{matrix}
0 & K(x) \\
1/\rho(x) & 0
\end{matrix}\right] \hspace{0.2in} \textnormal{and} \hspace{0.2in}
q(x,t) = \left[\begin{matrix}
p(x,t) \\ u(x,t)
\end{matrix}\right],
\end{align}
gives us the equation $q_t(x,t) + A(x)q_x(x,t) = 0$. 
Here $\rho(x)$ is the density and $K(x)$ is the bulk modulus. The eigenvalues of
$A(x)$ are $\pm c(x)$ where $c(x) = \sqrt{K(x)/\rho(x)}$ is the speed of sound.
We also define the impedance $Z(x) = \sqrt{K(x)\rho(x)} = \rho(x)c(x)$.  The
eigenvectors of $A(x)$ depend on the impedance; see the Appendix for more
details and solution to the Riemann problem for this problem and its adjoint.

\subsection{Example 1: Constant Impedance in One Dimension}
Let $t_0 = 0$, $t_f = 34$, $a = -12$, $b = 12$,
\begin{align*}
\rho (x) = \left\{
     \begin{array}{lr}
       1 & \hspace{0.3in}\textnormal{if } x < 0,\\
       4 & \textnormal{if } x > 0,
     \end{array}
   \right.
\qquad \text{and} \qquad
   K (x) = \left\{
     \begin{array}{lr}
       4 & \hspace{0.3in}\textnormal{if } x < 0,\\
       1 &\textnormal{if } x > 0,
     \end{array}
   \right.
\end{align*}
so that 
\begin{align*}
   c (x) = \left\{
     \begin{array}{lr}
       2 & \hspace{0.3in}\textnormal{if } x < 0,\\
       0.5 &\textnormal{if } x > 0,
     \end{array}
   \right.
\qquad \text{and} \qquad
Z(x) \equiv 2.
\end{align*}
Since the impedance is constant across $x=0$, so are the eigenvectors.
The sound speed changes and the waves are deformed as they pass
through the interface, but there will be no reflected waves. 
As initial data for $q(x,t)$ we take two wave packets in pressure,  
one on each side of the interface, 
and zero velocity everywhere. The wave packets in pressure are given by
\begin{equation}\label{eq:1d_qic}
p(x,0) = e^{-\beta_r(x-3)^2} \sin(f_rx) + e^{-\beta_l(x+2.5)^2} \sin(f_l x)
\end{equation}
with $\beta_l = 20$, $\beta_r = 5$, $f_l = 20$ and $f_r = 3$. 
As time progresses, the wave packets split into 
equal left-going and right-going waves which interact with the walls and 
the interface giving reflected and transmitted waves (recall that there will be no 
reflected waves off of the interface 
because the impedance is constant throughout the domain). 

Suppose that 
we are interested in the accurate estimation of the pressure near some point
$x_p$.  Then we can use the functional
 $J = \int_{a}^{b}\alpha \exp \left(-\hat{\beta} \left(x -
x_p\right)^2\right) p(x,t_f)\,dx$, which corresponds to \cref{eq:J_general} with
\begin{align}
\varphi (x) = \left[ \begin{matrix}
\alpha \exp \left(-\hat{\beta} \left(x - x_p\right)^2\right)\\ 0
\end{matrix}
\right]. \label{eq:phi_1d}
\end{align}
For this example we take $x_p = 7.5$, $\hat{\beta} = 50$,  and 
\begin{align}
\alpha = \sqrt{\hat{\beta}/\pi}
\end{align}
to normalize the Gaussian so it has mass 1 and represents an averaging of the
pressure in a small region around $x_p$.
\ignore{
Define
\begin{align*}
\hat{q}(x,t_f) = \left[ \begin{matrix}
\hat{p}(x,t_f) \\ \hat{u}(x,t_f)
\end{matrix}\right] = \varphi(x),
\end{align*}
and note that \cref{eq:q_integral} holds for this problem.
}
If we define the adjoint solution by solving
\begin{align}\label{eq:1d_adjoint}
\hat{q}_{t} + \left(A^T(x)\hat{q}\right)_{x} &= 0 &&x \in [a,b],
\hspace{0.1in}t_f \geq t  \geq t_0\\
\hat{u}(a,t) = 0, \hspace{0.1in} \hat{u}(b,t) &= 0 &&t_f \geq t  \geq t_0, \nonumber
\end{align}
then the second and third terms in \cref{eq:intbyparts} vanish and
we are left with \cref{eq:q_equality},
which is the expression that allows us to use the inner product of the
adjoint and forward problems at each time step to determine what regions will
influence the point of interest at the final time. 

\begin{figure}[h!]
\begin{minipage}[b]{0.45\linewidth}
\includegraphics[width=\textwidth]{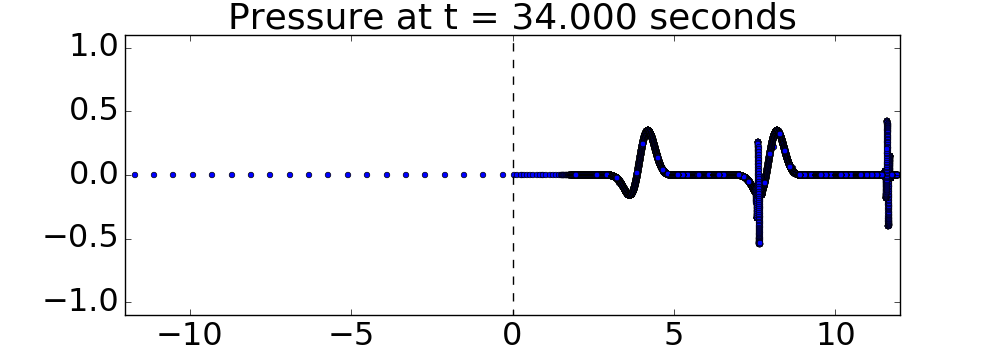}
\end{minipage}
\vspace{0.2cm}
\begin{minipage}[b]{0.45\linewidth}
\includegraphics[width=\textwidth]{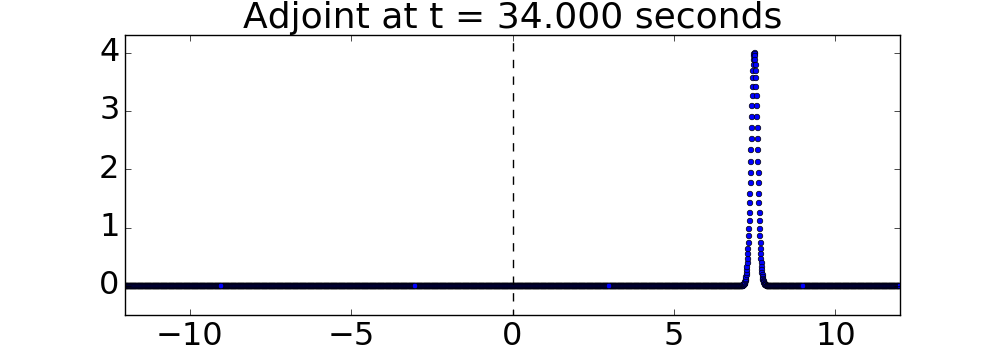}
\end{minipage}
\begin{minipage}[b]{0.45\linewidth}
\includegraphics[width=\textwidth]{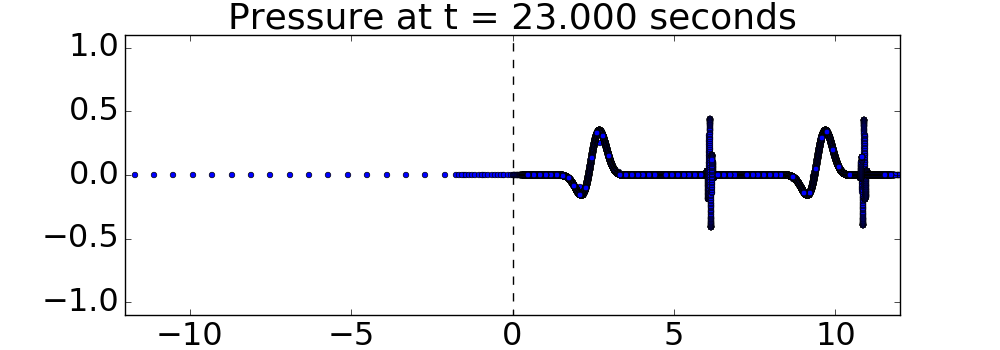}
\end{minipage}
\vspace{0.2cm}
\begin{minipage}[b]{0.45\linewidth}
\includegraphics[width=\textwidth]{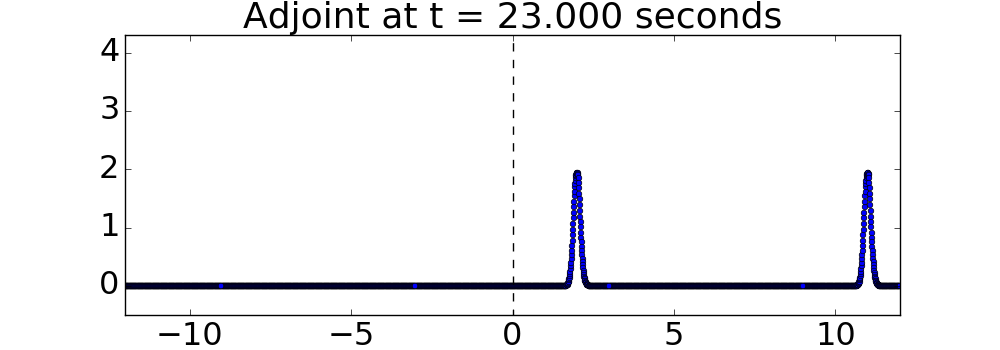}
\end{minipage}
\begin{minipage}[b]{0.45\linewidth}
\includegraphics[width=\textwidth]{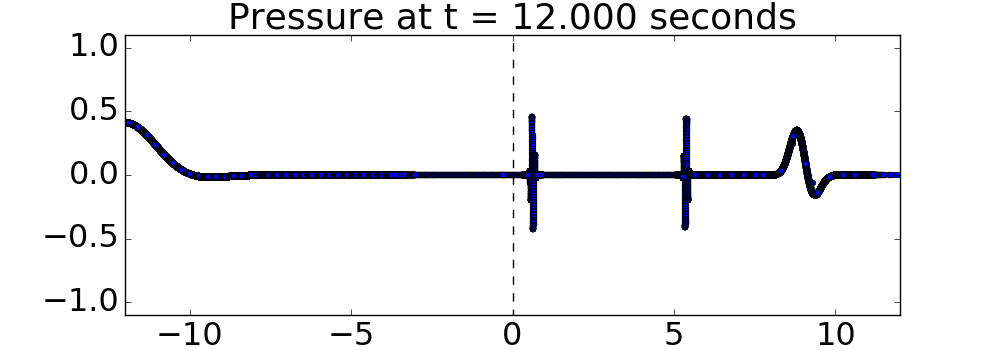}
\end{minipage}
\begin{minipage}[b]{0.45\linewidth}
\includegraphics[width=\textwidth]{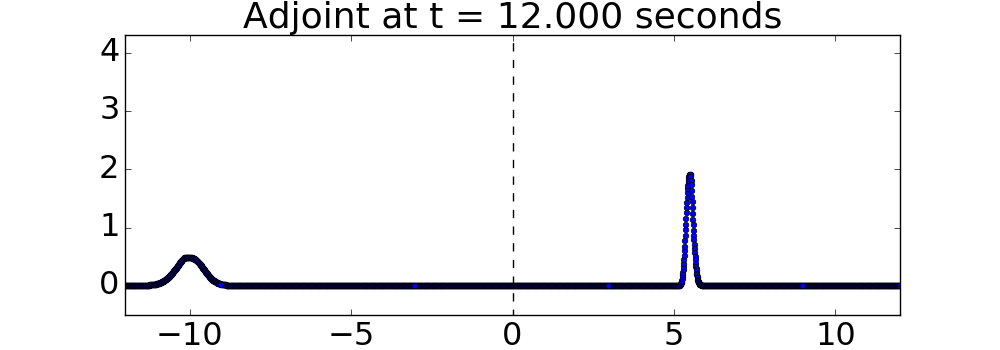}
\end{minipage}
\begin{minipage}[b]{0.45\linewidth}
\includegraphics[width=\textwidth]{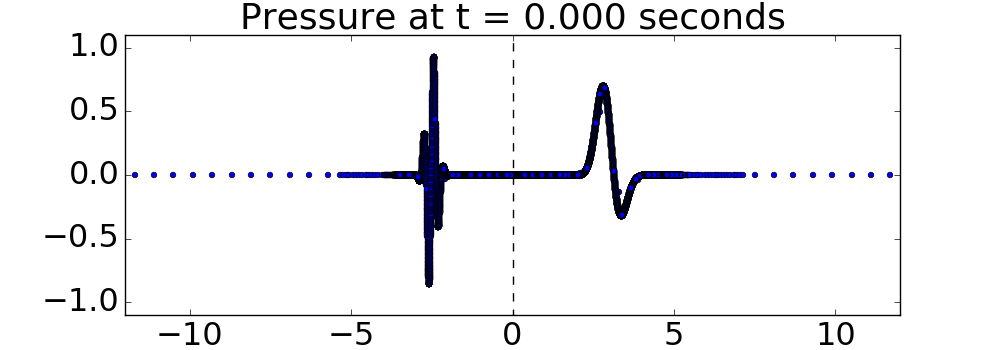}
\end{minipage}
\begin{minipage}[b]{0.45\linewidth}
\includegraphics[width=\textwidth]{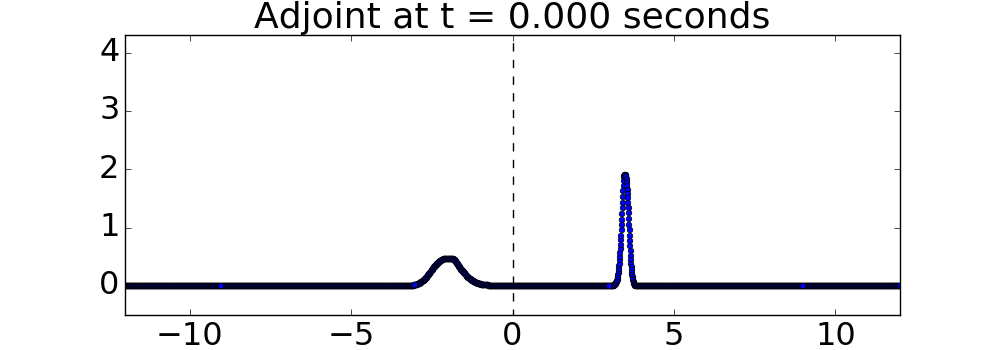}
\end{minipage}
 \caption{%
 Two wave packets in the forward problem (on the left, with initial data at
$t=0$, in the bottom plot)
 and a Gaussian hump in the 
 adjoint problem (on the right, with ``initial'' data at $t=34$, in the top plot),
 interacting with an interface at $x = 0$ (dashed line). There is a 
 change of sound speed at the interface, but the impedance is constant 
 throughout the domain. 
 The solutions are plotted as points, whose density varies due to AMR in the
 left plots.}
\label{fig:ex1_timeseries}
\end{figure}

As the ``initial'' data for this problem we set $\hat{q}(x,t_f) = \varphi
(x)$.
As time progresses backwards, the Gaussian hump splits into equal left-going and right-going
waves which interact with the walls and the interface giving reflected
 waves off of the walls and transmitted waves through the interface.
 \Cref{fig:ex1_timeseries} shows the AMRClaw  
results for the pressure of both the forward and adjoint solutions at
four different times. Each data point represents the pressure 
(the $y$ axis) at the center of a grid point in space (the $x$ axis). 
The density of data points in $x$ is indicative of the resolution 
of the level of refinement in that location. For instance, in the 
top plot on the left side of  \cref{fig:ex1_timeseries}, the widely 
spaced data points on the left side of the domain indicate that 
the region is covered by a very coarse grid while the tightly spaced 
data points on the right side of the domain indicate that the 
region is covered by a very fine grid. 
Both the forward and adjoint problems are run using $t_f = 34$, so the forward problem is
run from $t= 0$ to $t = 34$ and the adjoint problem is run from $t= 34$ to $t
= 0$.

\begin{figure}[h!]
  \centering
    \includegraphics[width=0.7\textwidth]{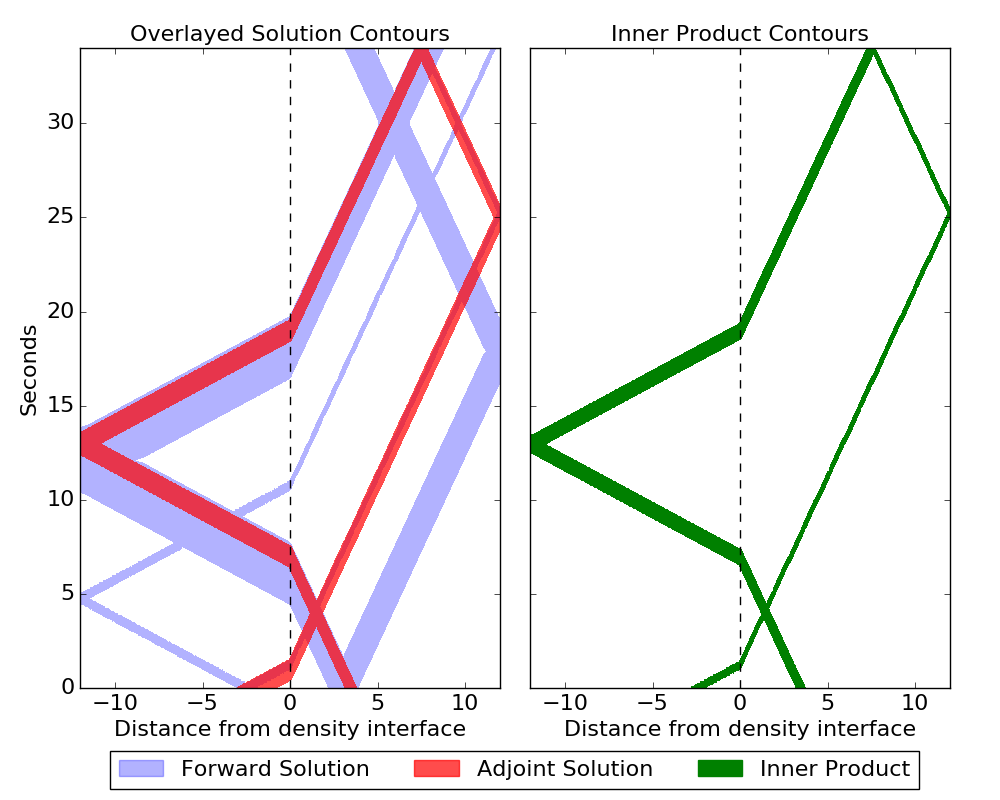}
    \caption{Space-time plots showing the forward and adjoint solutions
    from Example 1
    overlaid in the left plot.  Colored regions are where the magnitude of the
    first component of $q(x,t)$ or $\hat{q}(x,t)$ is above some threshold.
    The plot on the right shows regions where $|\hat{q}^T(x,t)q(x,t)|$ is
    above some tolerance,
    indicating where the forward problem should be
    refined if we are using adjoint-magnitude flagging. 
    The time axis is the same for both plots.}
    \label{fig:ex1_norms}
\end{figure}

To better visualize how the waves are moving through the domain, it is helpful
to look at the data in the $x$-$t$ plane as shown in \cref{fig:ex1_norms}. 
For \cref{fig:ex1_norms}, the
horizontal axis is the position, $x$, and the vertical axis is time. The left plot 
shows in blue the locations where $1-$norm of $q(x,t)$ is greater than or
equal to $10^{-2}$ and in red the locations where the $1-$norm
of $\hat{q}(x,t)$ is greater than or equal to $10^{-2}$. The right plot 
shows in green the locations where the inner product $\hat{q}^T(x,t)q(x,t)$ 
is greater than or equal to $10^{-2}$. This indicates which portions of the wave 
in the forward problem will actually reach our point of interest, and are 
exactly the regions that will be refined when using the adjoint-magnitude 
flagging method.

\subsection{Example 2: Variable Impedance in One Dimension}
As a slightly more complex example, 
let $t_0 = 0$, $t_f = 52$, $a = -12$, $b = 12$,
\begin{align*}
\rho (x) = \left\{
     \begin{array}{lr}
       1 & \hspace{0.3in}\textnormal{if } x < 0,\\
       2 & \textnormal{if } x > 0,
     \end{array}
   \right.
\qquad \text{and} \qquad
   K (x) = \left\{
     \begin{array}{lr}
       4 & \hspace{0.3in}\textnormal{if } x < 0,\\
       1 & \textnormal{if } x > 0,
     \end{array}
   \right.
\end{align*}
so that
\begin{align*}
   c (x) = \left\{
     \begin{array}{lr}
       2 & \hspace{0.3in}\textnormal{if } x < 0\\
       0.5 & \textnormal{if } x > 0
     \end{array}
   \right.
\qquad \text{and} \qquad
   Z (x) = \left\{
     \begin{array}{lr}
       4 & \hspace{0.3in}\textnormal{if } x < 0,\\
       0.5 & \textnormal{if } x > 0.
     \end{array}
   \right.
\end{align*}
In this case there is a jump in impedance at $x=0$, 
which will result in there being both transmitted and reflected waves 
at the interface. 
As initial data for $q(x,t)$ we take two wave packets in pressure 
one on each side of the interface, 
and zero velocity. The wave packets in pressure are the same as 
in example 1, and given by \cref{eq:1d_qic}. 
As time progresses, the wave packets split into 
equal left-going and right-going waves which interact with the walls and 
the interface giving reflected and transmitted waves. 
In contrast to example 1, many more waves arise in the solution as time
evolves due to the generation of new waves at each reflection.

\begin{figure}[h!]
  \centering
    \includegraphics[width=0.7\textwidth]{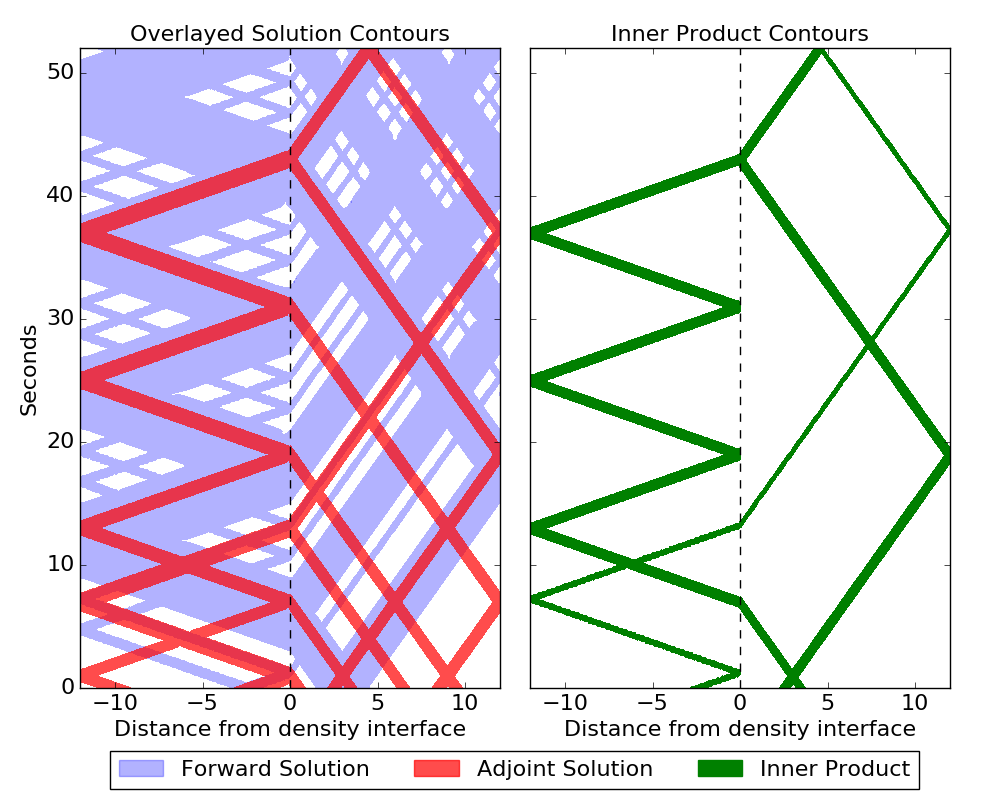}
    \caption{Space-time plots showing the forward and adjoint solutions
    from Example 2
    overlaid in the left plot.  Colored regions are where the magnitude of the
    first component of $q(x,t)$ or $\hat{q}(x,t)$ is above some threshold.
    The plot on the right shows regions where $|\hat{q}^T(x,t)q(x,t)|$ is
    above some tolerance,
    indicating where the forward problem should be
    refined if we are using adjoint-magnitude flagging. 
    The time axis is the same for both plots.}
    \label{fig:ex2_norms}
\end{figure}

For this example we suppose that 
we are interested in the accurate estimation of the pressure around
$x_p=4.5$ at $t=52$, chosen so that several waves in the forward problem
must be resolved to accurately 
capture the solution (while many other waves do not need
to be resolved; see \cref{fig:ex2_norms}).  
As ``initial'' data for the adjoint problem, $\hat{q}(x,t_f) = \varphi (x)$,
we can use the same functional \cref{eq:phi_1d} as in
Example 1, with this value of $x_p$, and with
$\alpha$, and $\hat\beta$ defined as in the previous example. 
As time progresses backwards in the adjoint solution, 
the hump splits into equal left-going and right-going
waves that interact with the walls and the interface giving both reflected
and transmitted waves.
Both the forward and adjoint problems are run using $t_f = 52$, so the forward problem is
run from $t= 0$ to $t = 52$ and the adjoint problem is run from $t= 52$ to $t
= 0$.

We will again visualize how the waves are moving through the domain by 
looking at the data in the $x$-$t$ plane as shown in \cref{fig:ex2_norms}. 
For \cref{fig:ex2_norms}, the
horizontal axis is the position, $x$, and the vertical axis is time. The left plot 
shows in blue the locations where $1-$norm of $q(x,t)$ is greater than or
equal to $10^{-5}$ and in red the locations where the $1-$norm
of $\hat{q}(x,t)$ is greater than or equal to $10^{-3}$. The right plot 
shows in green the locations where the inner product $\hat{q}^T(x,t)q(x,t)$ 
is greater than or equal to $10^{-3}$. This indicates which portions of the wave 
in the forward problem will actually reach our point of interest, and are 
exactly the regions that will be refined when using the adjoint-magnitude 
flagging method.

\begin{figure}[h!]
\begin{minipage}[b]{0.45\linewidth}
\includegraphics[width=\textwidth]{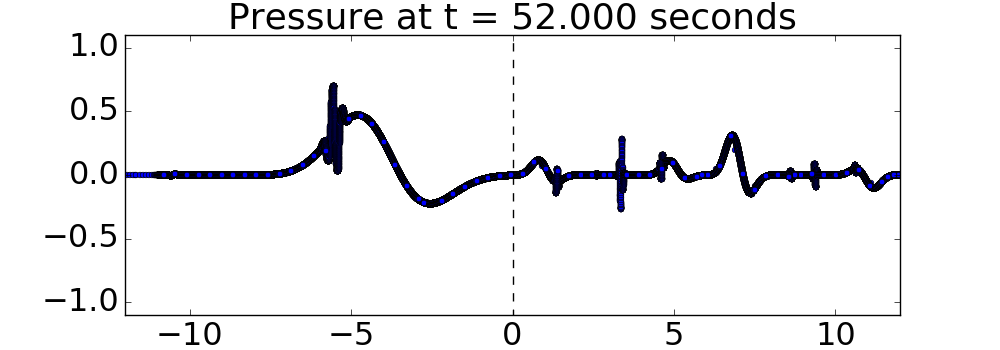}
\end{minipage}
\vspace{0.2cm}
\begin{minipage}[b]{0.45\linewidth}
\includegraphics[width=\textwidth]{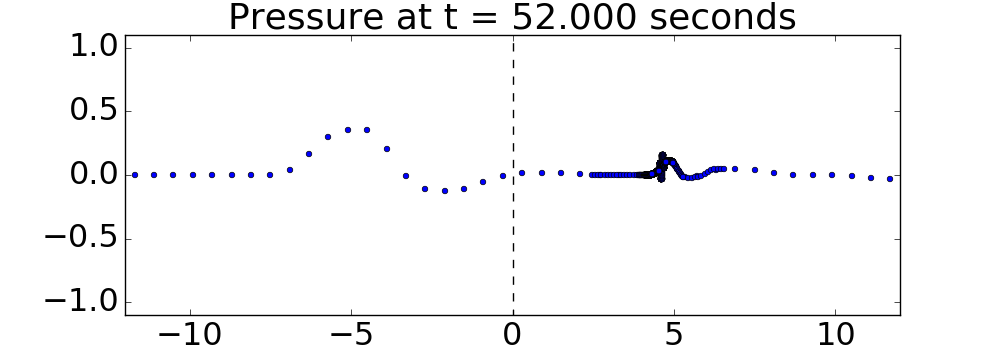}
\end{minipage}
 \caption{%
 The final solution computed in example 2 when using two different refinement
strategies.
 On the left: using difference-flagging. 
 On the right:  using adjoint-magnitude flagging, with a functional $J$ chosen to
resolve only the waves that affect the solution near $x_p=4.5$ at this time.}
\label{fig:ex2_finalsol_diff}
\end{figure}

\begin{figure}[h!]
\begin{minipage}[b]{0.45\linewidth}
\includegraphics[width=\textwidth]{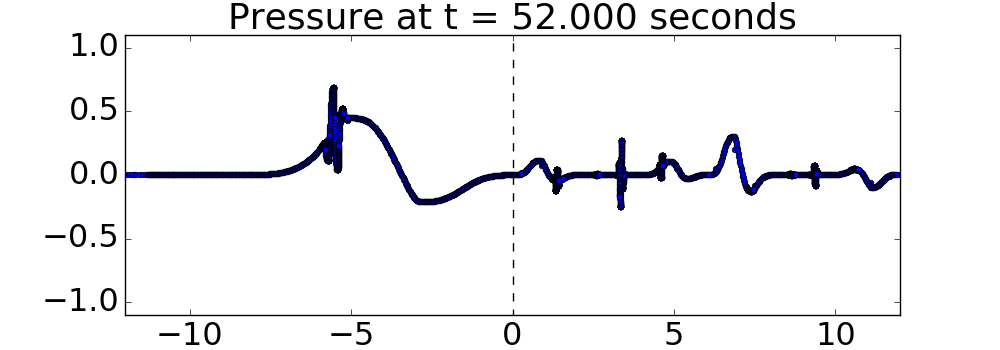}
\end{minipage}
\vspace{0.2cm}
\begin{minipage}[b]{0.45\linewidth}
\includegraphics[width=\textwidth]{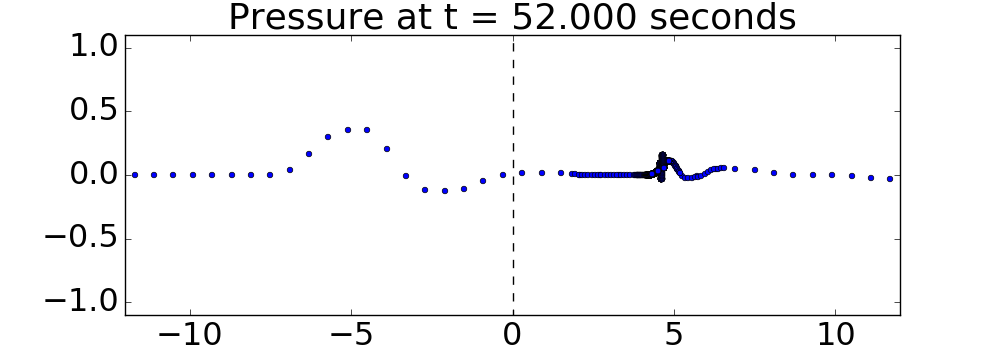}
\end{minipage}
 \caption{%
  The final solution computed in example 2 when using two different refinement
strategies.
 On the left: using error-flagging. 
 On the right:  using adjoint-error flagging, with a functional $J$ chosen to
resolve only the waves that affect the solution near $x_p=4.5$ at this time.
}
\label{fig:ex2_finalsol_err}
\end{figure}

To better visualize the impact these different flagging methods 
have on the final solution, \cref{fig:ex2_finalsol_diff} 
shows the final solution when using difference flagging and 
when using adjoint-magnitude flagging to achieve the same level of accuracy 
for our functional $J$. Similarly, \cref{fig:ex2_finalsol_err} 
shows the final solution when using error flagging and when 
using adjoint-error flagging to achieve the same level of accuracy in $J$. 
Note that the adjoint flagging methods have only 
focused on accurately capturing the waves at the 
location of interest, whereas the difference and error flagging 
techniques have captured all of the waves in the domain. 

\begin{figure}[h!]
\begin{minipage}[b]{0.95\linewidth}
\includegraphics[width=\textwidth]{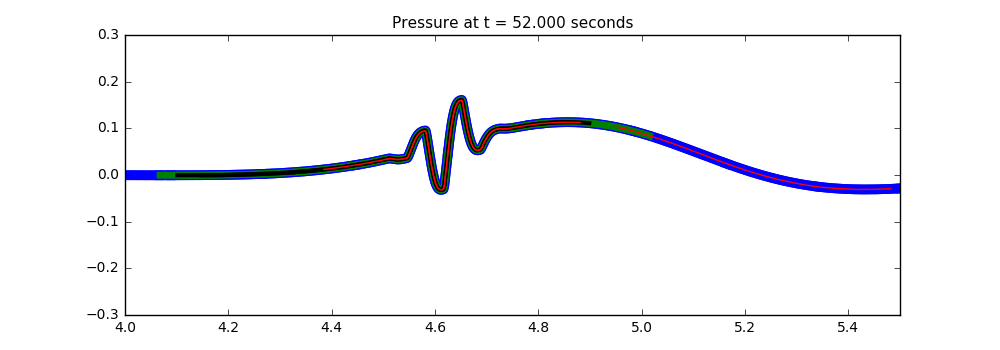}
\end{minipage}
 \caption{%
  The final solution computed in example 2 for all four refinement strategies, 
  in a zoomed in region of the domain. Difference-flagging 
  shown in blue, error-flagging shown in red, adjoint-magnitude flagging shown in 
  green, and adjoint-error flagging shown in black. Note that the lines are so 
  overlapped that they are almost impossible to tell apart. 
}
\label{fig:ex2_compare}
\end{figure}

As an aid to visualizing how the different methods capture the solution at the 
region of interest, \cref{fig:ex2_compare} shows the solution for the four different 
flagging methods for a zoomed-in area around the region of interest. Difference-flagging 
  is shown in blue, error-flagging is shown in red, adjoint-magnitude is flagging shown in 
  green, and adjoint-error flagging is shown in black. Note that although all four solutions 
  (from the four different flagging methods) are show, they are almost impossible to 
  distinguish from one another in the figure. In this figure, only the solution on refinement 
  levels 4 and 5 is shown. Note that the area of the domain that is resolved 
  to refinement levels 4 and 5 is smaller for the adjoint-flagging methods 
  than for their non-adjoint flagging counterparts, as expected. 
  
\subsection{Computational Performance}\label{sec:1d_preformance}
Recall that we are considering four different flagging methodologies: 
difference-flagging, error-flagging, 
adjoint-magnitude flagging, and adjoint-error flagging. To compare the 
results from these methods we will take into 
account 
\begin{itemize}
\item the placement of AMR patches throughout the domain,  
\item the amount of CPU time required, 
\item the number of cell updates required, and 
\item the accuracy of the computed solution
\end{itemize}
for each flagging method. 

Let us begin by considering the placement of AMR levels throughout the domain. 
Since the adjoint-flagging methods take into account only the portions of the wave 
that will impact the region of interest during our time range of interest we expect that 
using adjoint-flagging will result in less of the domain being covered by fine 
resolution grids --- where we are making the assumption that some of the 
waves will not ultimately have an effect on our region of interest. 

In both of the examples we have presented 
this is in fact the case. Consider the contour plots shown in  
\cref{fig:ex1_norms} and \cref{fig:ex2_norms}. If we are using difference-flagging 
then anywhere the forward solution is large (all of the areas in blue) would be 
flagged. In contrast, if we are using adjoint-magnitude flagging all of the 
areas where the inner product $\hat{q}^T(x,t)q(x,t)$ is large (all of 
the areas in green) would be flagged. Therefore, more of the domain 
would be flagged when using difference-flagging, resulting in more of 
the domain being covered by finer levels of refinement. 
For both of these one-dimensional examples a total of five levels of refinement 
are used for the forward problem, 
starting with $40$ cells on the coarsest level and with a refinement 
ratio of $6$ from each level to the next. 
So, the finest level in the forward problem 
corresponds to a fine grid with 51,840 cells. 
Where these finer levels of refinement are placed, of course, varies 
based on the flagging method being used.
The adjoint problem was solved on a relatively 
coarse grid with $3000$ cells, and 
no adaptive mesh refinement. 

\begin{figure}[ht]
 \centering
  \subfigure[Levels generated by using difference-flagging.]{
  \includegraphics[width=0.475\textwidth]{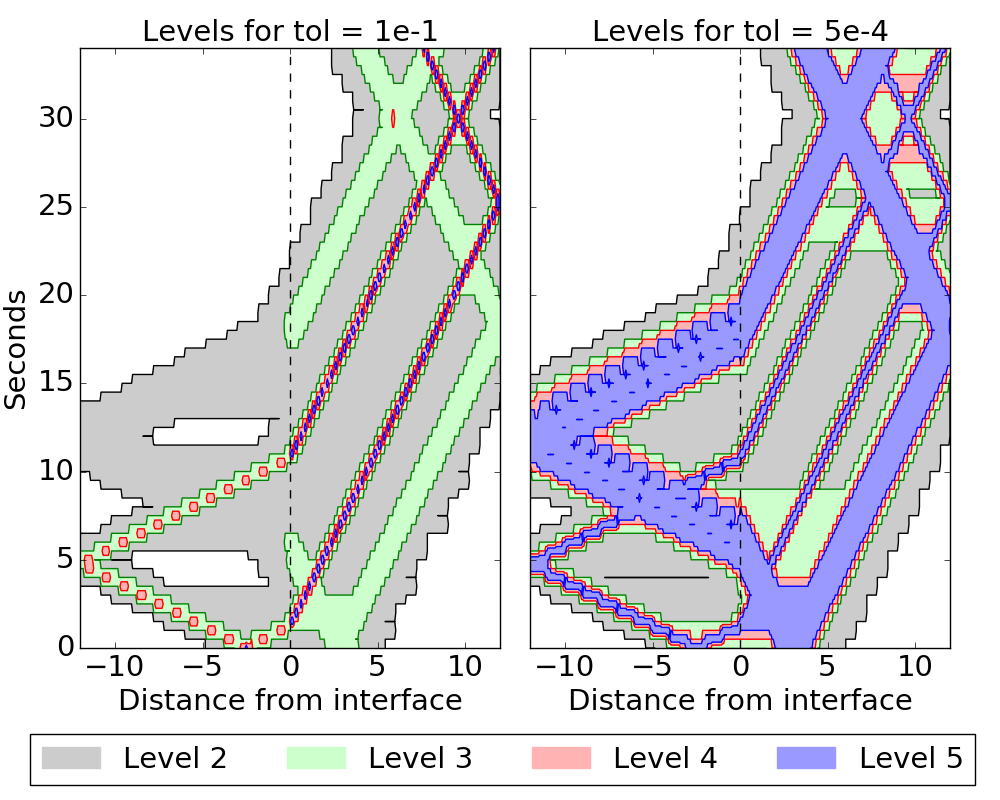}
   }
  \subfigure[Levels generated by using adjoint-magnitude flagging.]{
  \includegraphics[width=0.475\textwidth]{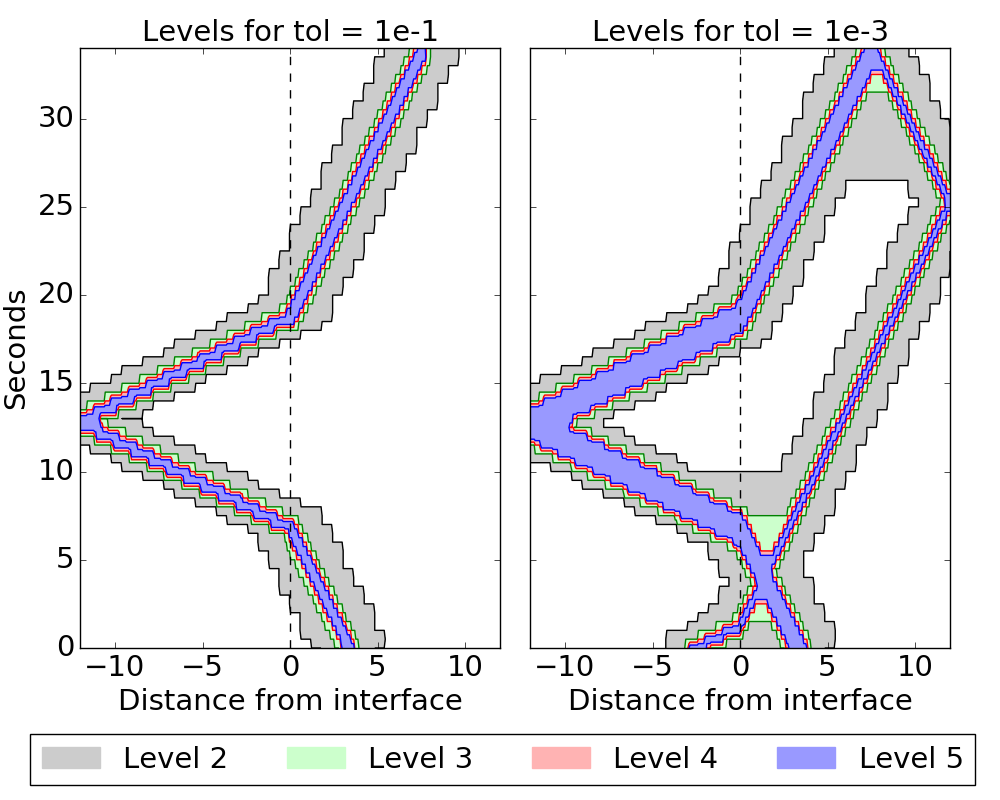}
   }
 \caption[Optional caption for list of figures]{%
   Levels used for the difference-flagging and adjoint-magnitude flagging methods for example 1. 
  Tolerances used for the two difference-flagging runs shown are $10^{-1}$ and $5 \times 10^{-4}$. 
  Tolerances used for the two adjoint-magnitude  
  flagging runs shown are $10^{-1}$ and $10^{-3}$.   
  A larger level number means greater resolution: so 
  level 5 for example has a finer resolution than level 4. 
  Recall that level 1 covers the entire domain, and is therefore not shown.}
   \label{fig:ex1_mag_levels}
\end{figure}

\begin{figure}[ht]
 \centering
  \subfigure[Levels generated by using error-flagging.]{
  \includegraphics[width=0.475\textwidth]{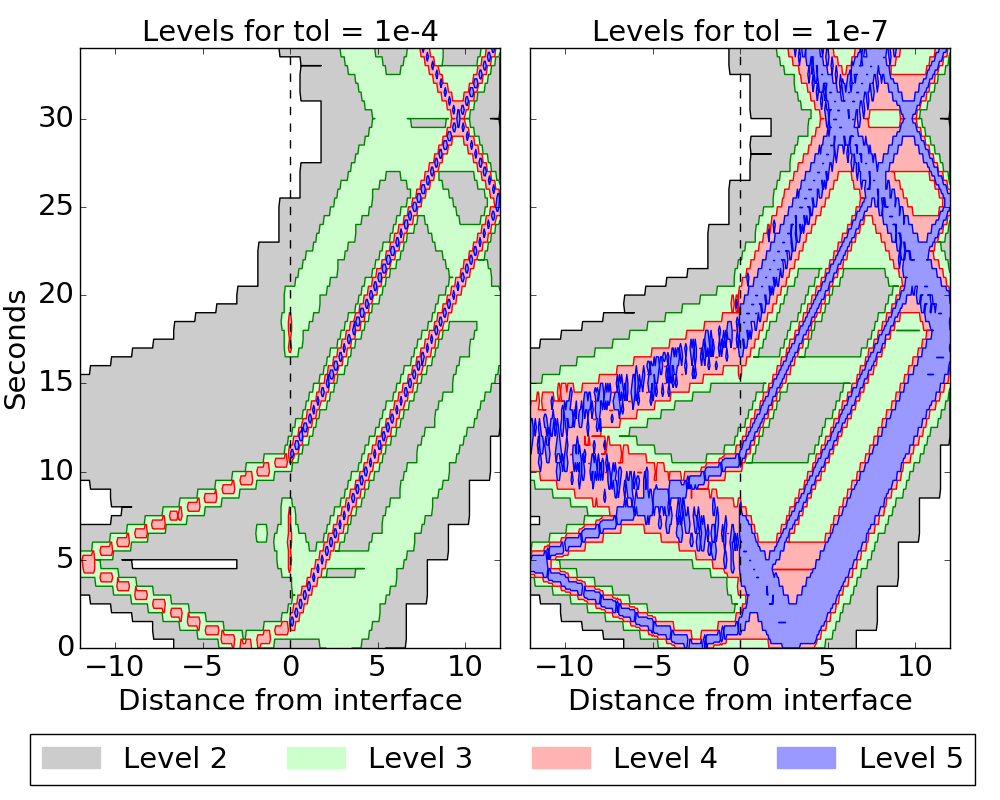}
   }
  \subfigure[Levels generated by using adjoint-error flagging.]{
  \includegraphics[width=0.475\textwidth]{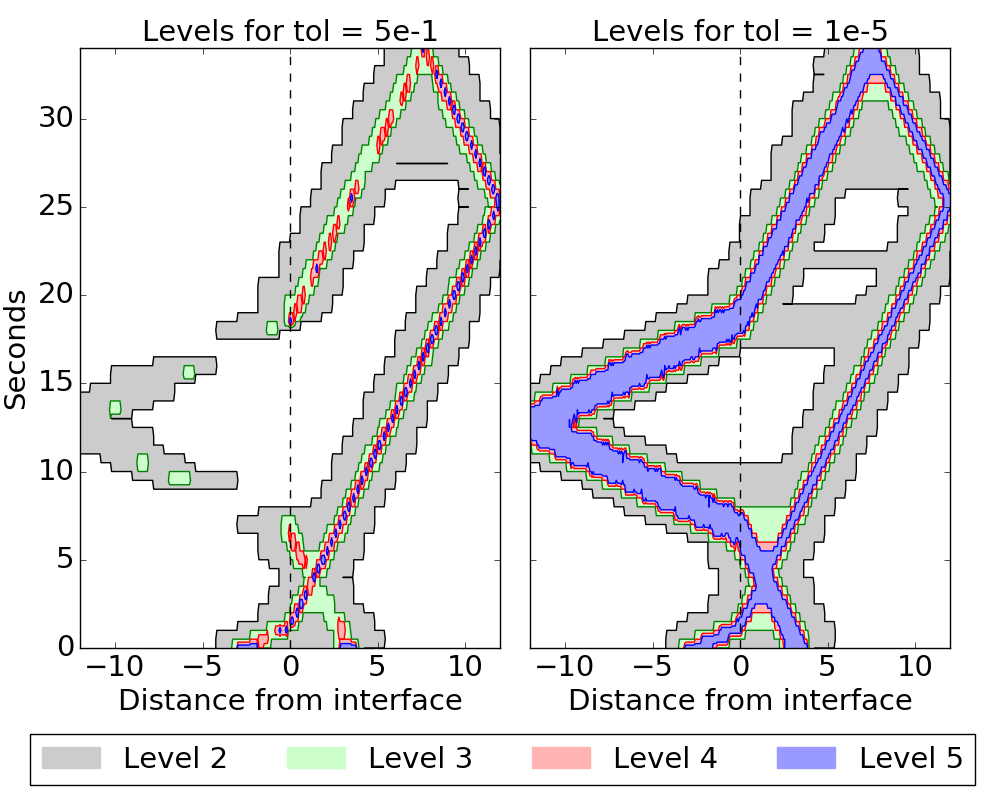}
   }
 \caption[Optional caption for list of figures]{%
  Levels used for the error-flagging and adjoint-error flagging methods for example 1. 
  Tolerances used for the two error-flagging runs shown are $10^{-4}$ and $10^{-7}$. 
  Tolerances used for the two adjoint-error 
  flagging runs shown are $5 \times 10^{-1}$ and $10^{-5}$. 
  Recall that level 1 covers the entire domain, and is therefore not shown.}
   \label{fig:ex1_error_levels}
\end{figure}

\Cref{fig:ex1_mag_levels} shows the levels of refinement for 
difference-flagging (on the left) and adjoint-magnitude flagging (on the right) 
for two different tolerances
in example 1. The larger tolerance used for both of these flagging methods 
is $10^{-1}$ which yields 
a final error in the functional of interest of about $10^{-2}$ for both methods.
The smaller tolerance used for difference flagging is $5 \times 10^{-4}$ 
and for adjoint-magnitude flagging is $10^{-3}$, both of 
which yield an error of about 
$10^{-6}$ in the functional of interest.
Note that these two flagging methods 
both consider the value of the pressure in the cell when determining 
whether or not the cell should be flagged --- difference-flagging flags 
the cell based on the relationship between the pressure in the current 
cell and the pressure in the adjacent cells, and adjoint-magnitude flagging 
flags the cell based on the value of the inner product between $q(x,t)$ 
and $\hat{q}(x,t)$ in the cell. As expected, when the adjoint problem 
is taken into account the result is a smaller region of the domain 
begin refined to achieve a comparable level of error.

Similarly, \cref{fig:ex1_error_levels} shows the levels of 
refinement for error-flagging (on the left) and adjoint-error 
flagging (on the right) for two different tolerances in example 1. 
The larger tolerance used for error-flagging 
is $10^{-4}$ and for adjoint-error flagging is $5 \times 10^{-1}$, both of which yield 
a final error  in the functional of interest of about $10^{-2}$.
The smaller tolerance used for error-flagging 
is $10^{-7}$ and for adjoint-error flagging is $10^{-5}$, both of which yield 
an error of about $10^{-6}$ in the functional of interest. 
These two flagging methods consider the estimated error in the solution 
in each cell when determining whether or not the cell should be 
flagged. However, the adjoint-error flagging method takes the 
additional step of 
considering the inner product of this error estimate
with the adjoint solution. Again, when the adjoint problem 
is taken into account the result is a smaller region of the domain 
being refined to achieve a comparable final error.

The tolerances shown were selected so that 
these figures would allow for easy comparison between the four methods. 
All of the larger tolerances chosen yield an error of about $10^{-2}$ 
in our functional of interest for each of the methods, and all of the 
smaller tolerances chosen yield an error of about $10^{-6}$ in our functional
 of interest. Therefore, it is appropriate to compare the extent of the
 refined regions present in these figures. Note that even though 
 the same accuracy is achieved, the placement of 
 regions of greater refinement varies between these four methods. 

\begin{figure}[ht]
 \centering
  \subfigure[Levels generated by using difference-flagging for example 2.]{
  \includegraphics[width=0.475\textwidth]{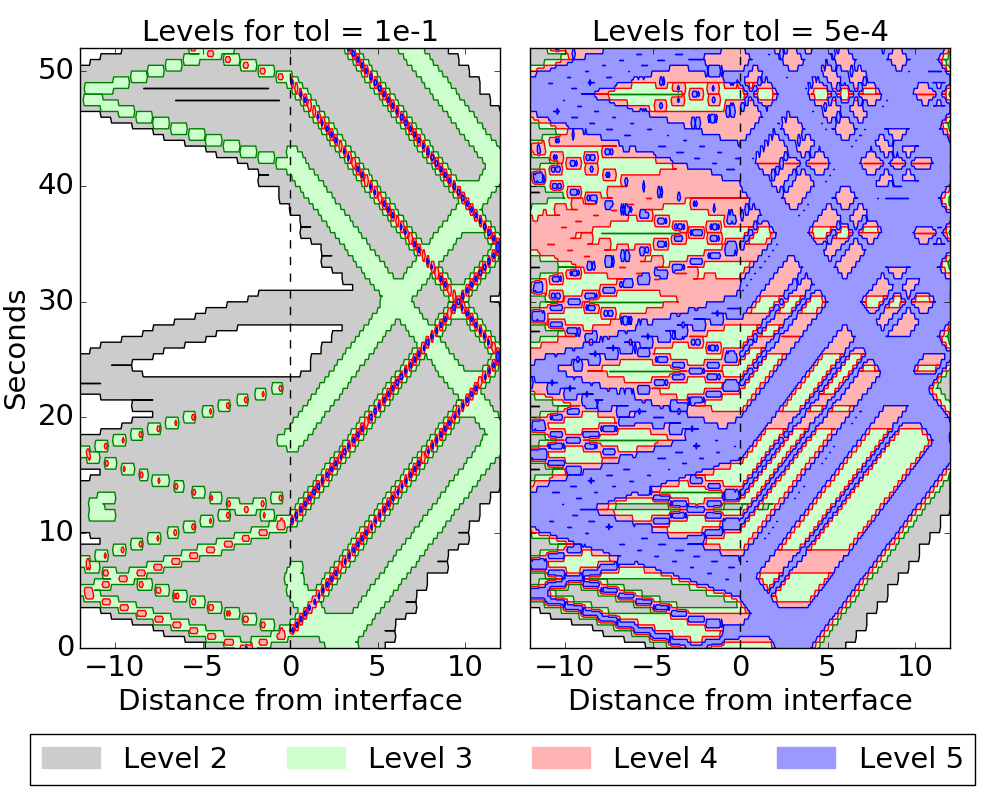}
   }
  \subfigure[Levels generated by using adjoint-magnitude flagging for example 2.]{
  \includegraphics[width=0.475\textwidth]{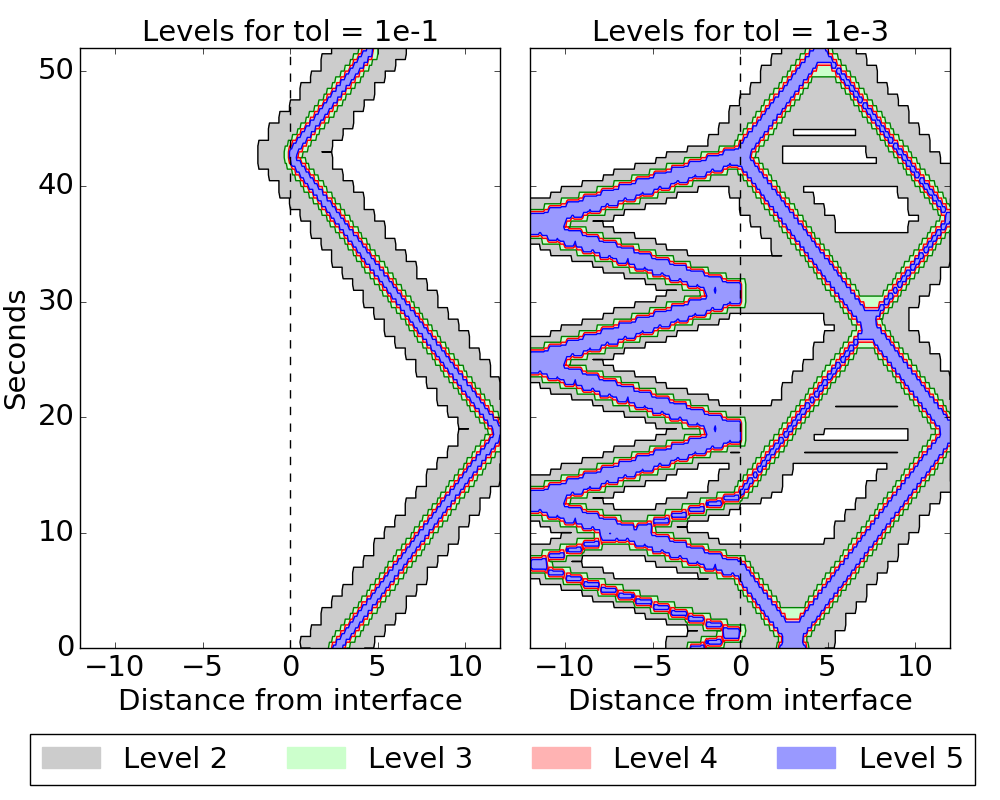}
   }
 \caption[Optional caption for list of figures]{%
   Levels used for the difference-flagging and adjoint-magnitude flagging methods for example 2. 
  Tolerances used for the two difference-flagging runs shown are $10^{-1}$ and $5 \times 10^{-4}$. 
  Tolerances used for the two adjoint-magnitude  
  flagging runs shown are $10^{-1}$ and $10^{-3}$.   
  Recall that level 1 covers the entire domain, and is therefore not shown.}
   \label{fig:ex2_mag_levels}
\end{figure}

\begin{figure}[ht]
 \centering
  \subfigure[Levels generated by using error-flagging for example 2.]{
  \includegraphics[width=0.475\textwidth]{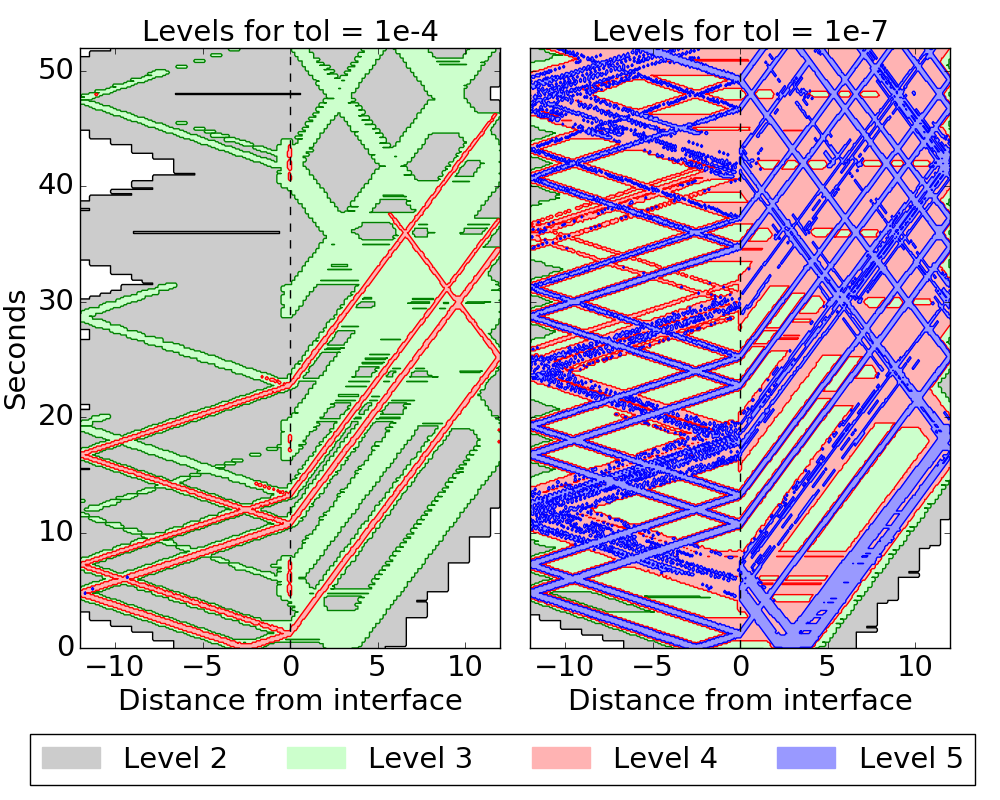}
   }
  \subfigure[Levels generated by using adjoint-error flagging for example 2.]{
  \includegraphics[width=0.475\textwidth]{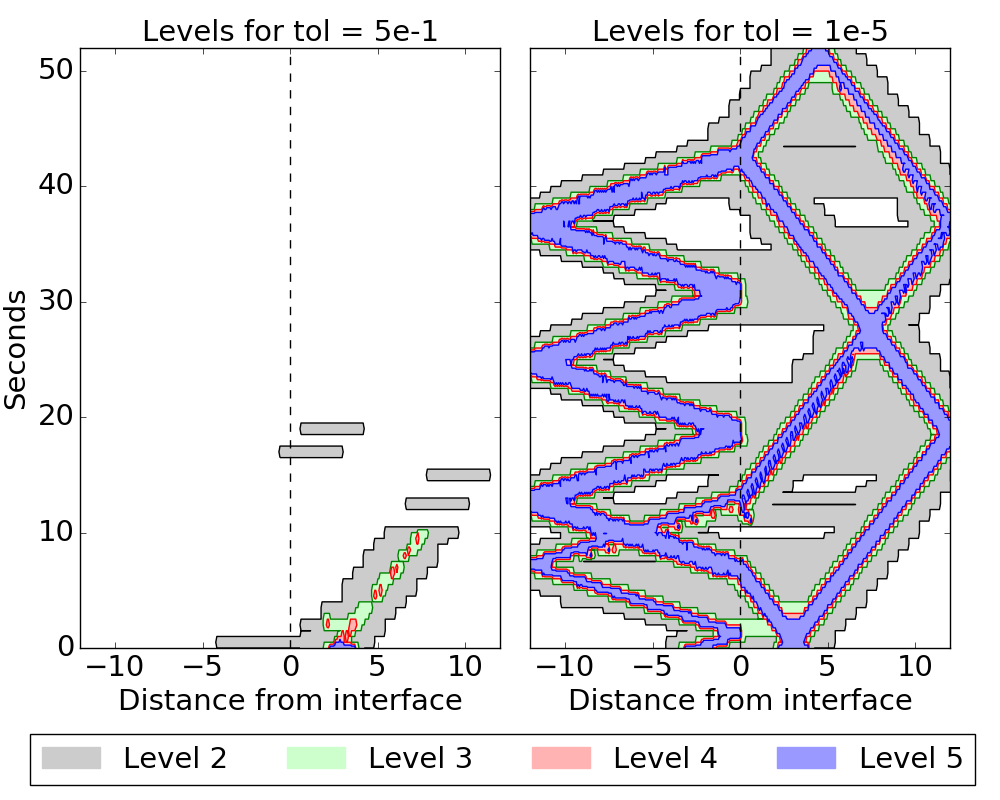}
   }
 \caption[Optional caption for list of figures]{%
  Levels used for the error-flagging and adjoint-error flagging methods for example 2. 
  Tolerances used for the two error-flagging runs shown are $10^{-4}$ and $10^{-7}$. 
  Tolerances used for the two adjoint-error 
  flagging runs shown are $5 \times 10^{-1}$ and $10^{-5}$. 
  Recall that level 1 covers the entire domain, and is therefore not shown.}
   \label{fig:ex2_error_levels}
\end{figure}

For example 2, \cref{fig:ex2_mag_levels} shows the levels of refinement for 
difference-flagging (on the left) and adjoint-magnitude flagging (on the right) 
using the same tolerances that were shown for example 1. 
\Cref{fig:ex2_error_levels} shows the levels of 
refinement for error-flagging (on the left) and adjoint-error 
flagging (on the right) for example 2, 
where we have chosen the same tolerances that were shown 
for example 1.
For this example we again see that when the adjoint problem 
is taken into account the result is a smaller region of the domain 
begin refined.

\begin{figure}[h!]
\begin{minipage}[b]{0.45\linewidth}
\includegraphics[width=\textwidth]{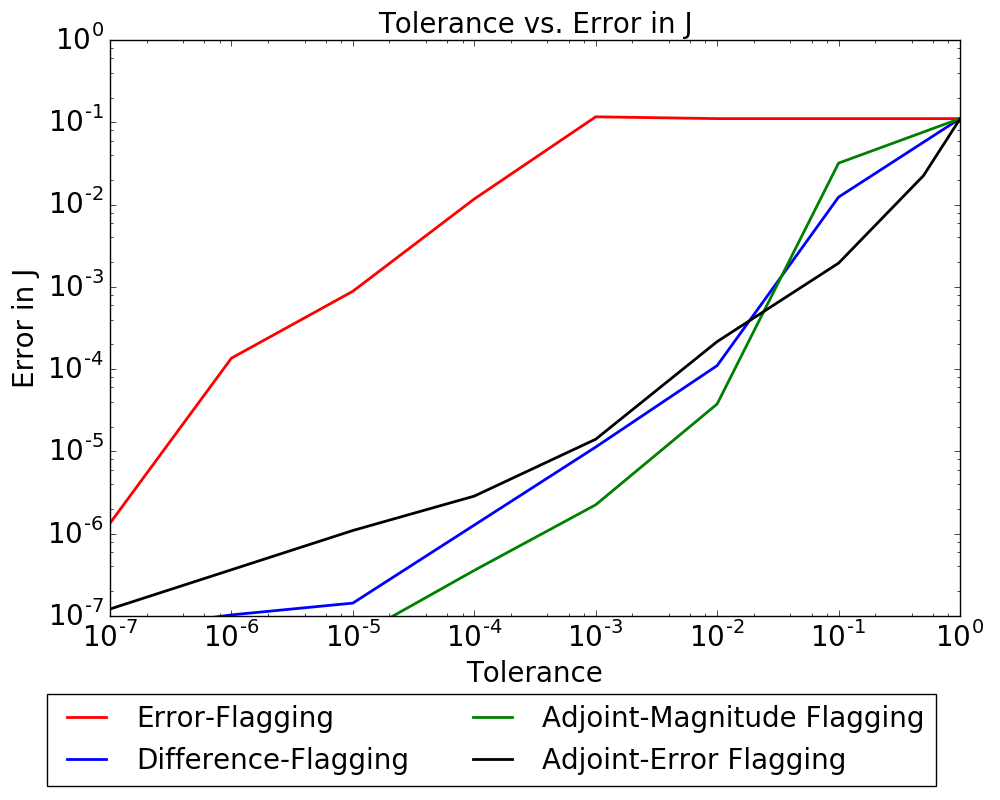}
\end{minipage}
 \caption{%
 Accuracy for the different flagging methods on 
 example 1. Shown is the error in the functional J for each tolerance value 
 for the various flagging methods.
 }
\label{fig:ex1_tolvserror}
\end{figure}

Also of interest is the amount of work that is required for 
each of the flagging methods we are considering, and 
the level of accuracy that resulted from that work. 
\Cref{fig:ex1_tolvserror} shows the error in our functional of 
interest given in \cref{eq:J} for example 1 
as the tolerance is varied for the four flagging 
methods we are considering.
Difference-flagging was used with a tolerance of $10^{-12}$ to 
compute a very fine grid solution to our forward problem for this example. 
This solution was used to compute a fine grid value of our functional 
of interest, $J_\text{fine}$. 
The error 
between the calculated value for $J$ and $J_\text{fine}$ is shown. 

Note that the magnitude of the error is 
nearly the same for a given tolerance when using 
either difference-flagging or adjoint-magnitude flagging although, 
as we will see, the amount  
of work required is significantly less when using adjoint-flagging. 
Also note that the magnitude of the error when using adjoint-error 
flagging is consistently less than the magnitude of the tolerance being used, 
as we expected based on \cref{sec:adjErrorFlag}. This 
means that adjoint-error flagging allows the user to enforce 
a certain level of accuracy on the final functional of interest by 
selecting a tolerance of the desired order. This capability has not 
existed in Clawpack previous to this work. 
Recall that the tolerance set by the user is being used 
to evaluate whether or not some given quantity is above that set 
tolerance. However, the quantity that is being evaluated is different 
for each of the flagging methods we are considering. Therefore, 
other than noting the general trends of each line
in this figure, one cannot reach conclusions on relative merits by comparing 
one line with the others in this plot.

\begin{figure}[h!]
\begin{minipage}[b]{0.45\linewidth}
\includegraphics[width=\textwidth]{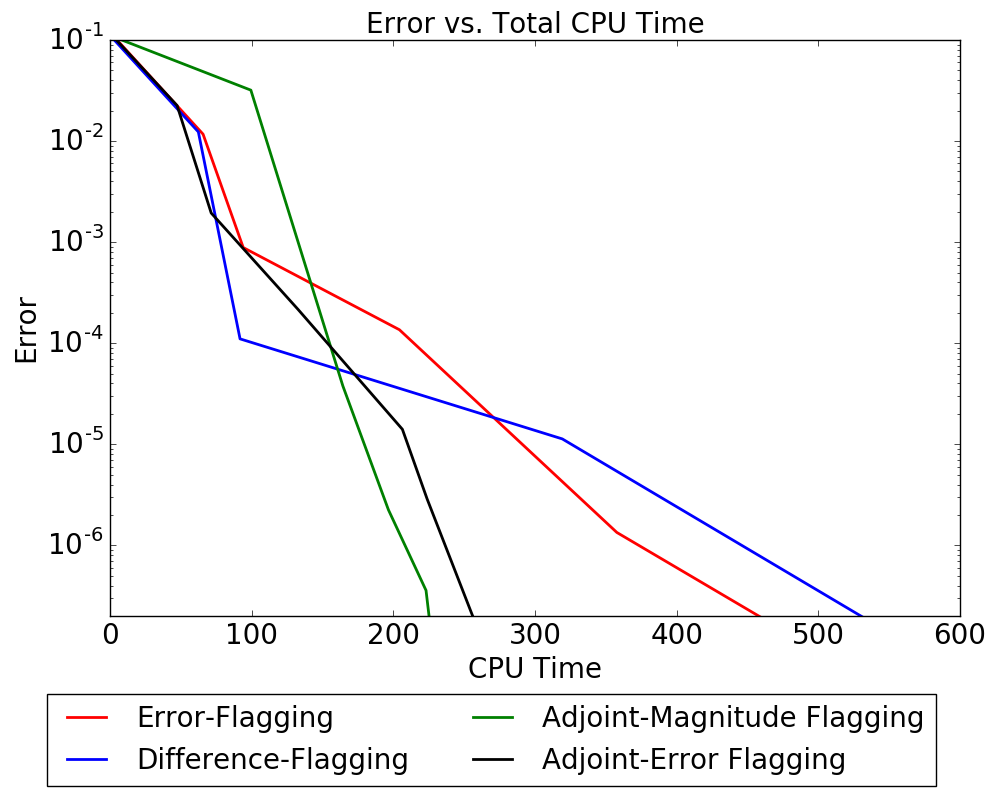}
\end{minipage}
\vspace{0.2cm}
\begin{minipage}[b]{0.45\linewidth}
\includegraphics[width=\textwidth]{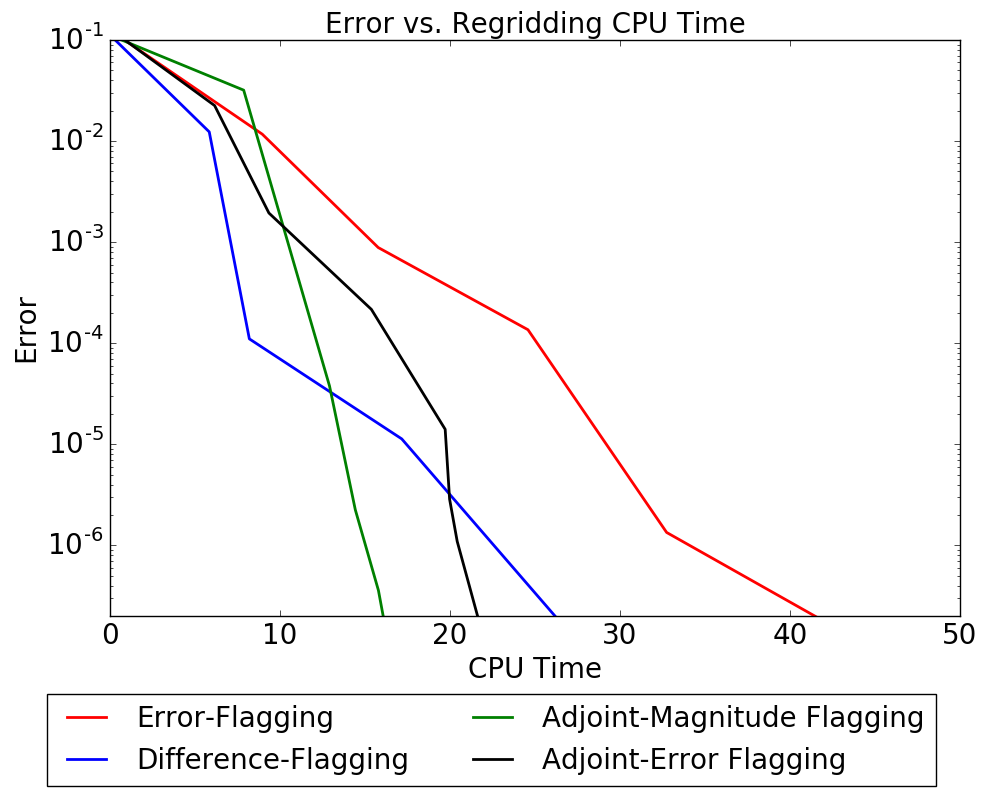}
\end{minipage}
 \caption{%
Performance measures for the different flagging methods on 
 example 1. 
 On the left: the total CPU time (in seconds) required vs. the 
 accuracy achieved.
 On the right: the regridding CPU time (in seconds) required vs. the 
 accuracy achieved. CPU times were found by averaging the CPU time 
 over fifteen runs.}
\label{fig:ex1_accuracy}
\end{figure}

\begin{figure}[h!]
\begin{minipage}[b]{0.45\linewidth}
\includegraphics[width=\textwidth]{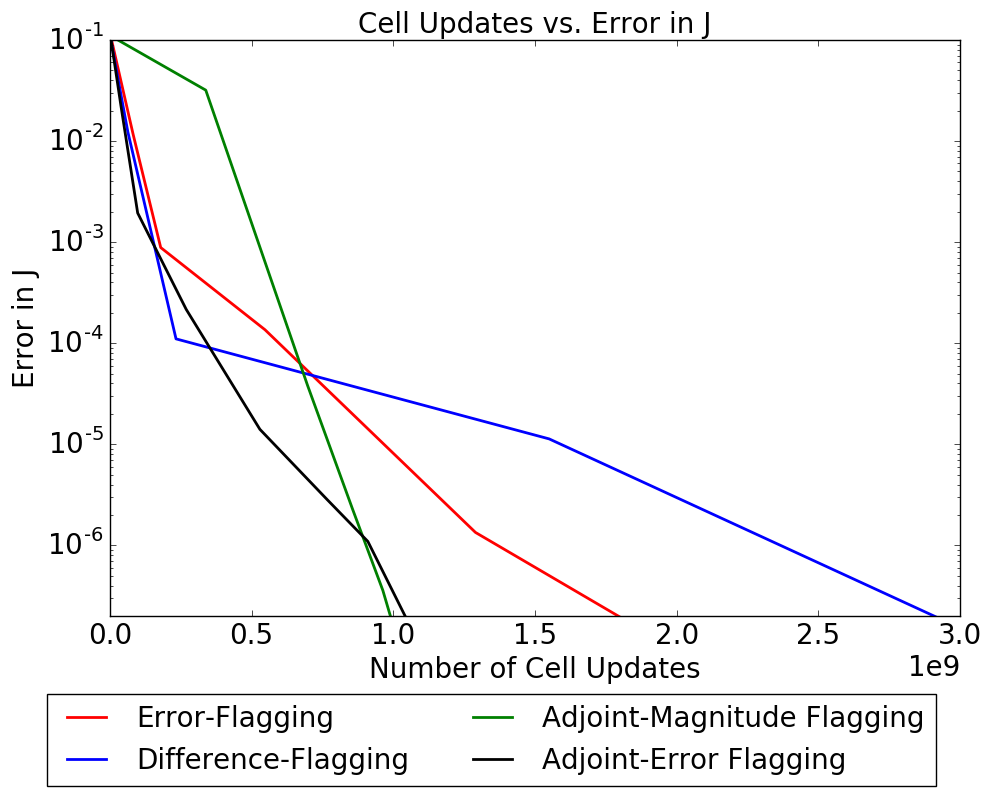}
\end{minipage}
 \caption{%
Number of cell updates calculated vs. the accuracy achieved 
for the four different flagging methods begin considered. Note 
that the number of cell updates axis is multiplied by $1$e$9$.
 }
\label{fig:ex1_cellsvserror}
\end{figure}

We are better able to compare the methods to one another if we
consider the amount of computational time that is required to 
achieve a certain level of accuracy with each flagging method. 
The above examples were run on a quad-core laptop, 
and the OpenMP option of AMRClaw was enabled, which 
allowed all four cores to be utilized. \Cref{fig:ex1_accuracy} 
shows two measures of the amount of CPU time that was 
required for each method vs. the accuracy achieved. On the left, 
we have shown the total amount of CPU time used by the 
computation. This includes stepping the solution forward in time 
by updating cell values, 
regridding at appropriate time intervals, outputting the results, 
and other various overhead requirements. For each flagging 
method and tolerance this example was run fifteen times, 
and the average of the CPU times for those runs was used 
in this plot.  
Note that while 
the adjoint-error and adjoint-difference flagging  
have higher CPU time requirements for lower accuracy, 
these two methods quickly show their strength by maintaining 
a low CPU time requirement while increasing the accuracy of the 
solution. In contrast, the CPU time requirements for 
difference-flagging and error-flagging quickly increase when increased 
accuracy is required. 
Also remember that the adjoint-flagging methods do have the 
additional time requirement of solving the adjoint 
problem, which is not shown in these figures. For this example, 
solving the adjoint required about 10 seconds of CPU 
time, which is small compared to the time spent on the 
forward problem. This CPU time was found by taking the average 
time required over ten simulations. 

The right of \cref{fig:ex1_accuracy} shows the amount of CPU time 
spent on the regridding process vs. the accuracy of the functional $J$. 
Note that this is really where the differences between the four methods 
lie: the adjoint-flagging methods goal is to reduce the number of cells 
that are flagged while maintaining the accuracy of the solution. While 
we have seen from the left side of \cref{fig:ex1_accuracy} that 
this has certainly been successful in reducing the overall time required 
for computing the solution, we might expect that the time spent in the 
actual regridding process might be longer for the adjoint-flagging methods. 
This is due to the fact that adjoint-flagging requires various extra steps 
when determining whether or not to flag a cell (recall, for example, that 
we need to calculate the inner product between the adjoint and 
either the forward solution or the estimated error in the forward solution). However, 
as we can see from this figure, the amount of time spent in the 
regridding process is not significantly greater for the adjoint-flagging 
methods. In fact, for higher levels of accuracy the adjoint-flagging 
methods required less time for the regridding process than their 
non-adjoint counterparts because there are fewer fine grids. 

Another advantage of the adjoint-flagging methods is the fact that they 
have less memory requirements than the alternative flagging methods. 
As we have already stated, since the adjoint method allows us to safely ignore regions of the 
domain, there are fewer fine grids as observed from the fact that many fewer 
cell updates were required. \Cref{fig:ex1_cellsvserror} shows 
the number of cell updates required vs. the accuracy achieved for each of the 
flagging methods being considered. Note that the number of cell updates 
is significantly less for the adjoint-flagging methods than for their non-adjoint 
counterparts. This reduces the memory requirements for the computation. 
For this relatively small example memory usage when utilizing the non-adjoint 
flagging method did not become an issue, but it can quickly become a 
constraint for larger computations (in particular for three dimensional problems).

\section{Two-Dimensional Variable Coefficient Acoustics}\label{sec:2d_acoustics_example}
In two dimensions the variable coefficient linear acoustics equations are 
\begin{equation}\label{acou2d}
\begin{split} 
p_t(x,y,t) + K(x,y) \left(u_x(x,y,t) + v(x,y,t)_y\right) &= 0, \\
\rho(x,y) u_t(x,y,t) + {p_x(x,y,t)} &= 0, \\
\rho(x,y) v_t(x,y,t) + {p_y(x,y,t)} &= 0,
\end{split}
\end{equation} 
in the domain $x \in [a,b],~y \in [\alpha,\beta],~t > t_0$.
Setting 
\begin{align*}
A(x,y) = \left[ \begin{matrix}
0 & K(x,y) & 0 \\
1/\rho(x,y)& 0 & 0 \\
0 & 0 & 0
\end{matrix}\right],\hspace{0.1in}
B(x,y) = \left[ \begin{matrix}
0 & 0 & K(x,y) \\
0 & 0 & 0 \\
1/\rho(x,y)& 0 & 0
\end{matrix}\right], \hspace{0.1in}
q(x,y,t) = \left[\begin{matrix}
p(x,y,t) \\ u(x,y,t) \\ v(x,y,t)
\end{matrix}\right],
\end{align*}
gives us the equation $q_t(x,y,t) + A(x,y)q_x(x,y,t) + B(x,y)q_y(x,y,t) = 0$. 

As an example,
consider $t_0 = 0$, $t_f = 21$, $a = -8$, $b = 8$, $\alpha = -1$, 
$\beta = 11$, 
\begin{align*}
\rho(x,y) \equiv 1, 
\qquad \text{and} \qquad
   K (x,y) = \left\{
     \begin{array}{lr}
       4 & \hspace{0.3in}\textnormal{if } x < 0,\\
       1 & \textnormal{if } x > 0,
     \end{array}
   \right.
\end{align*}
so that
\begin{align*}
   c (x,y) = \left\{
     \begin{array}{ll}
       2 & \hspace{0.3in}\textnormal{if } x < 0,\hspace{0.04in}\\
       1 &  \hspace{0.3in}\textnormal{if } x > 0,
     \end{array}
   \right.
\qquad \text{and} \qquad
   Z (x,y) = \left\{
     \begin{array}{ll}
       2 & \hspace{0.3in}\textnormal{if } x < 0,\hspace{0.04in}\\
       1 &  \hspace{0.3in}\textnormal{if } x > 0.
     \end{array}
   \right.
\end{align*}
We use wall boundary conditions on the top and both sides, 
\begin{align}
&u(a,y,t) = 0, \hspace{0.1in} u(b,y,t) = 0 &&t  \geq 0, \\
& v(x,\beta,t) = 0 &&t
  \geq 0, \label{eq:acoustics2d_wallsBC}
\end{align}
and 
extrapolation boundary conditions on the bottom at $y = -1$.
The impedance jump at
$x = 0$ will result in reflected and transmitted waves emanating from 
this interface. Also, we will have reflected waves from the walls along 
both sides and the top of our domain. 

As initial
data for $q(x,y,t)$ we take a smooth radially symmetric hump in pressure, and
a zero velocity in both $x$ and $y$ directions. The initial hump in pressure is given by
\begin{align}
p(x,y,0) = \left\{
     \begin{array}{lr}
       3 + \cos \left(\pi \left(r - 0.5\right)/w\right) 
       & \hspace{0.3in}\textnormal{if } \left| r - 0.3 \right| \leq w,\\
       0 & \textnormal{otherwise.} \hspace{0.27in}
     \end{array}
   \right. \label{eq:2d_acoustics_p}
\end{align}
with $w = 0.15$ and $r = \sqrt{\left(x-0.5\right)^2 +\left( y - 1\right)^2}$. 
As time progresses, this hump in pressure will radiate outward symmetrically.

We will consider two different examples, one where there are just a few waves that 
intersect at our region of interest during our time range of interest, and a 
second where a larger number of waves contribute to our region and time 
range of interest. 

\subsection{Example 3: Capturing a few intersecting waves}
\begin{figure}[!htbp]
 \centering
  \subfigure[$t_f -$ 1 second]{
  \includegraphics[width=0.225\textwidth]{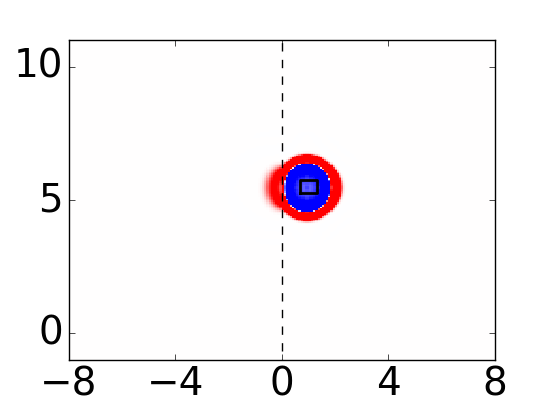}
   }
   \subfigure[$t_f -$ 6 seconds]{
  \includegraphics[width=0.225\textwidth]{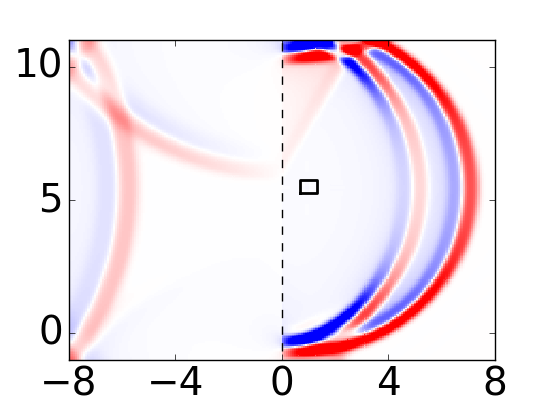}
   }
  \subfigure[$t_f -$ 15 seconds]{
  \includegraphics[width=0.225\textwidth]{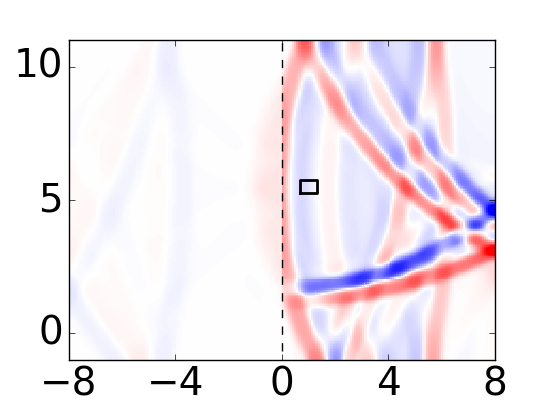}
   }
  \subfigure[$t_f -$ 18.5 seconds]{
  \includegraphics[width=0.225\textwidth]{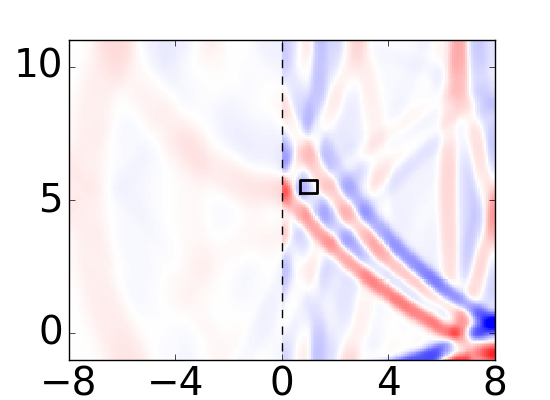}
   }
 \caption{%
   Computed results for two-dimensional acoustics adjoint problem 
   for example 3. Times shown are
   the number of seconds before the final time, since the ``initial'' conditions are
   given at the final time.
  The small rectangle in this and later figures shows the region of
interest defined by the functional $J$, and is the region where $I(x,y)=1$.
  The color scale goes from blue to red, and ranges between $-0.005$ and $0.005$.}
   \label{fig:2d_adjoint_ex3}
\end{figure}

Suppose that we are interested in the accurate estimation of the 
pressure in the area defined by a rectangle centered about 
$(x,y) = (1.0,5.5)$ at the final time $t_f = 21$.
The small rectangle in \Cref{fig:2d_adjoint_ex3} and later figures shows this 
region of interest.
Setting 
\begin{align*}
J = \int_{0.68}^{1.32}\int_{5.26}^{5.74}p(x,y,t_f)dy\,dx,
\end{align*}
the problem then requires that
\begin{align}
\varphi (x,y) = \left[ \begin{matrix}
 I(x,y) \\ 0 \\ 0
\end{matrix}
\right], \label{eq:phi_2d_ex3}
\end{align}
where
\begin{align}
I (x,y) = \left\{
     \begin{array}{ll}
       1 & \hspace{0.3in}\textnormal{if } 0.68 \leq x \leq 1.32
        \textnormal{ and } 5.26 \leq y \leq 5.74,\\
       0 &\hspace{0.3in} \textnormal{otherwise.}\hspace{1.61in}
     \end{array}
   \right. \label{eq:delta_2d_ex3}
\end{align}
Define
\begin{align*}
\hat{q}(x,y,t_f) = \left[ \begin{matrix}
\hat{p}(x,y,t_f) \\ \hat{u}(x,y,t_f) \\ \hat{v}(x,y,t_f)
\end{matrix}\right] = \varphi(x,y),
\end{align*}
and note that for any time $t_0 < t_f$ we have
\begin{align*}
\int_{t_0}^{t_f}\int_a^b\int_{\alpha}^{\beta}  \hat{q}^T\left( q_t + A(x,y)q_x +
B(x,y)q_y\right) dy\,dx\,dt = 0.
\end{align*}
Integrating by parts yields the equation
\begin{align}
\left. \int_a^b  \int_{\alpha}^{\beta}  \hat{q}^Tq\,dy\,dx \right|^{t_f}_{t_0}
&+ \left. \int_{t_0}^{t_f}\int_{\alpha}^{\beta}  \hat{q}^TA(x,y)q\,dy\,dt \right|^{b}_{a}
+ \left. \int_{t_0}^{t_f}\int_{a}^{b}  
\hat{q}^TB(x,y)q\,dx\,dt \right|^{\beta}_{\alpha}\nonumber \\
- &\int_{t_0}^{t_f}\int_a^b\int_{\alpha}^{\beta}  q^T\left(\hat{q}_{t} +
\left(A^T(x,y)\hat{q}\right)_{x} + \left(B^T(x,y)\hat{q}\right)_{y}\right)
dy\,dx\,dt = 0.\label{eq:acousticseqn_2d_ex3}
\end{align}
Note that if we can define an adjoint problem such that all but the first term 
in this equation vanishes then we are left with
\begin{align*}
\int_a^b \int_{\alpha}^{\beta}  \hat{q}^T(x, y, t_f)q(x, y, t_f) dy\,dx
= \int_a^b \int_{\alpha}^{\beta} \hat{q}^T(x, y, t_0)q(x, y, t_0)dy\,dx,
\end{align*}
which is the expression that allows us to use the inner product of the
adjoint and forward problems at each time step to determine what regions will
influence the point of interest at the final time. 
We can accomplish this by  
defining the adjoint problem
\begin{align}
&\hat{q}_{t} + \left(A^T(x,y)\hat{q}\right)_{x} + \left(B^T(x,y)
\hat{q}\right)_{y} = 0
&&x \in [a,b], y \in [\alpha, \beta], \hspace{0.1in} t_f \geq t \geq t_0 \label{eq:acoustics2d_adjointwalls_ex3}
\end{align}
with the same boundary conditions as the forward problem. 

As the initial data for the adjoint $\hat{q}(x,y,t_f) = \varphi (x,y)$ 
we have a square pulse in pressure, which was described in 
\cref{eq:phi_2d_ex3} and \cref{eq:delta_2d_ex3}.
As time progresses backwards, waves radiate outward and reflect off the
walls as well as transmitting and reflecting off of the interface at $x = 0$. 
The adjoint problem was solved on a $200 \times 200$ grid 
without mesh refinement, and several snapshots of the 
solution are shown in \cref{fig:2d_adjoint_ex3}. 
The location of interest is outlined with a black box in the plots. 

Although we have focused on the functional $J$, and are considering 
the accuracy and error in our estimates of this functional, in 
actuality we are typically interested in the accurate estimation 
of the pressure in the forward solution at a particular spatial 
and temporal area. The selection of the initial data for the adjoint 
plays a role in determining the functional $J$, and dictates the 
relationship between the solution of the forward problem 
and our functional of interest. Therefore, the initial data we select for the 
adjoint problem is undoubtedly significant. For this example 
we have selected initial conditions for the adjoint problem such that 
our functional $J$ is simply the integral of the forward solution 
over our region of interest. For the one-dimensional examples shown 
earlier, a Gaussian initial condition for the adjoint was used. 
An exploration of how these initial conditions affect the accuracy 
of the forward solution is not conducted in this work, but 
is an interesting area for future work.

\begin{figure}[h!]
\begin{minipage}[c]{0.06\linewidth}
\hspace{0.1cm}
\end{minipage}
\begin{minipage}[c]{0.23\linewidth}
\centering
t = 2.5 seconds
\end{minipage}
\begin{minipage}[c]{0.23\linewidth}
\centering
t = 6 seconds
\end{minipage}
\begin{minipage}[c]{0.23\linewidth}
\centering
t = 15 seconds
\end{minipage}
\begin{minipage}[c]{0.23\linewidth}
\centering
t = 21 seconds
\end{minipage}\\
\begin{minipage}[c]{0.06\linewidth}
\begin{sideways}
\parbox{2.9cm}{\centering Fine Grid \\ Solution}
\end{sideways}
\end{minipage}
\begin{minipage}[c]{0.23\linewidth}
\includegraphics[width=\textwidth]{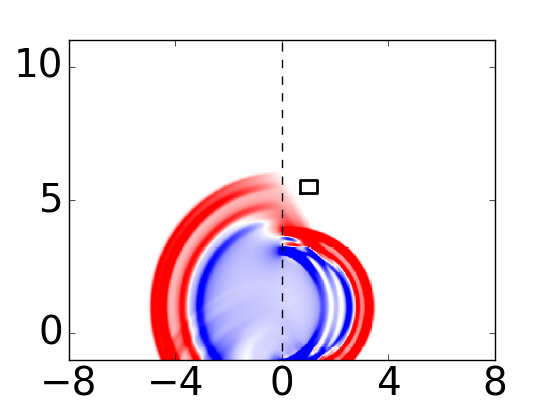}
\end{minipage}
\begin{minipage}[c]{0.23\linewidth}
\includegraphics[width=\textwidth]{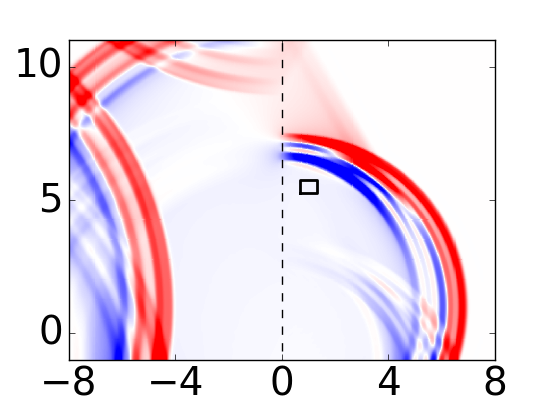}
\end{minipage}
\begin{minipage}[c]{0.23\linewidth}
\includegraphics[width=\textwidth]{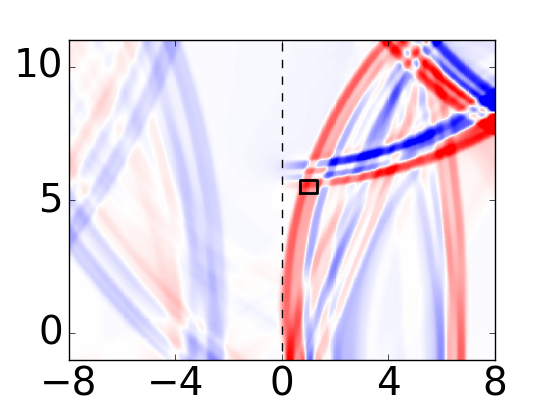}
\end{minipage}
\begin{minipage}[c]{0.23\linewidth}
\includegraphics[width=\textwidth]{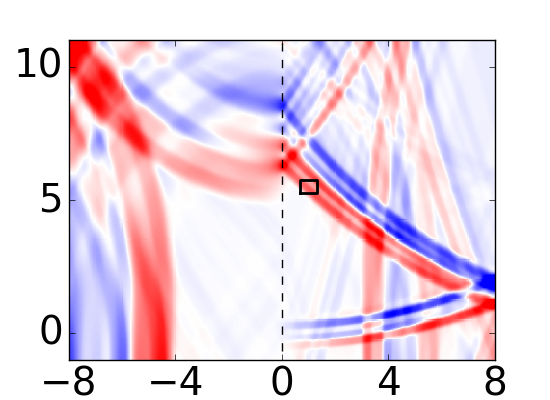}
\end{minipage}\\
\begin{minipage}[c]{0.06\linewidth}
\begin{sideways}
\parbox{2.9cm}{\centering Difference-Flagging \\ Grids}
\end{sideways}
\end{minipage}
\begin{minipage}[c]{0.23\linewidth}
\includegraphics[width=\textwidth]{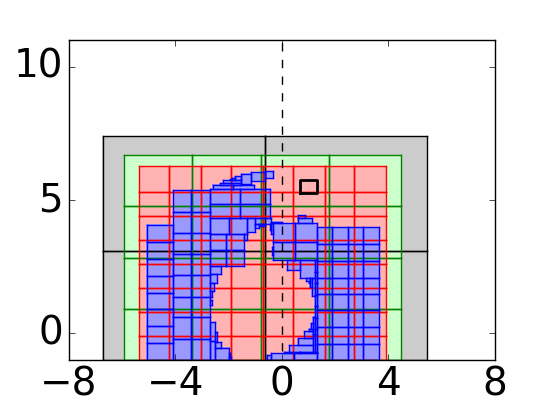}
\end{minipage}
\begin{minipage}[c]{0.23\linewidth}
\includegraphics[width=\textwidth]{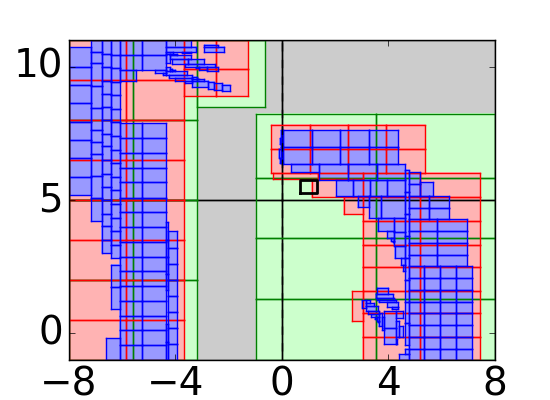}
\end{minipage}
\begin{minipage}[c]{0.23\linewidth}
\includegraphics[width=\textwidth]{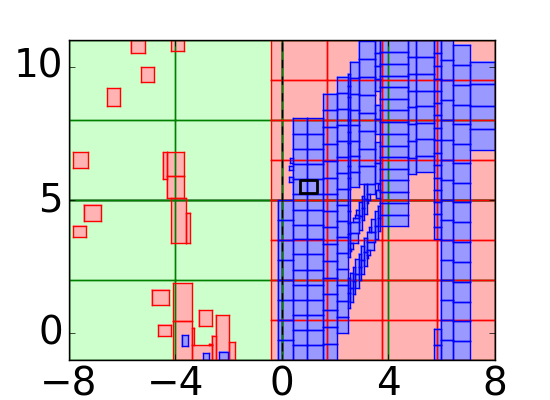}
\end{minipage}
\begin{minipage}[c]{0.23\linewidth}
\includegraphics[width=\textwidth]{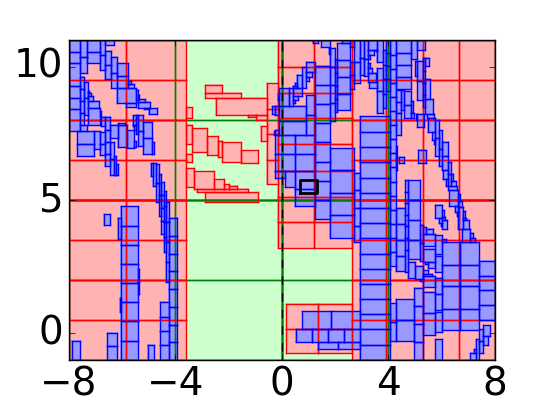}
\end{minipage}\\
\begin{minipage}[c]{0.06\linewidth}
\begin{sideways}
\parbox{2.9cm}{\centering Error-Flagging \\ Grids}
\end{sideways}
\end{minipage}
\begin{minipage}[c]{0.23\linewidth}
\includegraphics[width=\textwidth]{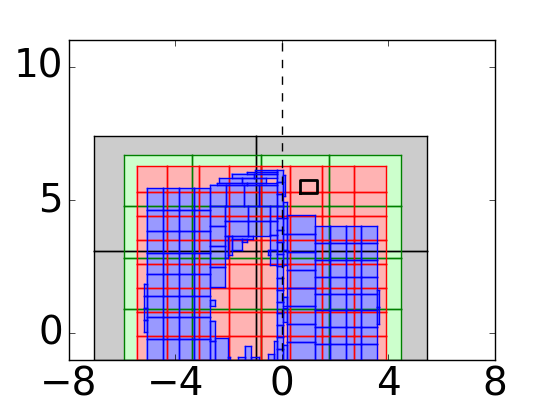}
\end{minipage}
\begin{minipage}[c]{0.23\linewidth}
\includegraphics[width=\textwidth]{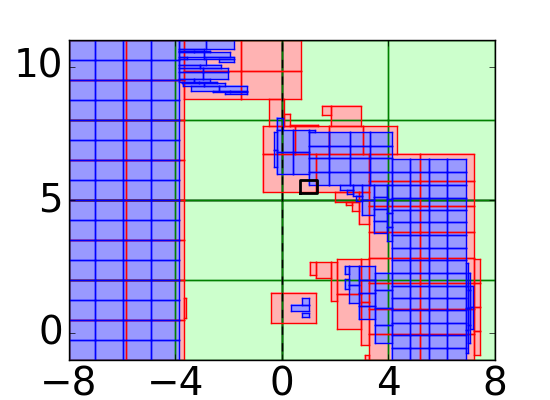}
\end{minipage}
\begin{minipage}[c]{0.23\linewidth}
\includegraphics[width=\textwidth]{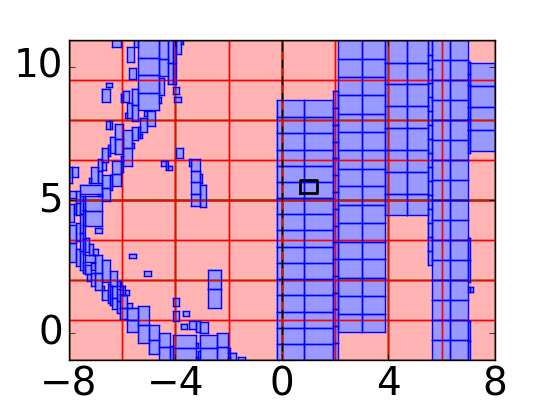}
\end{minipage}
\begin{minipage}[c]{0.23\linewidth}
\includegraphics[width=\textwidth]{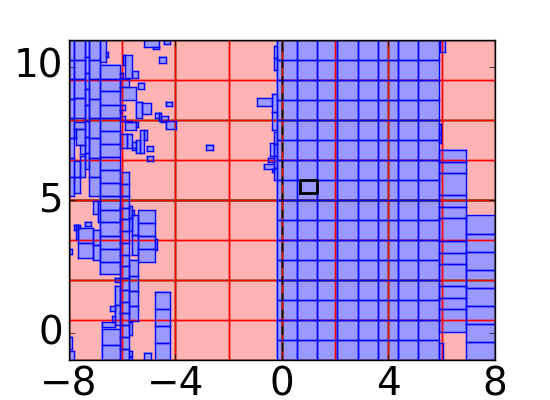}
\end{minipage}\\
\begin{minipage}[c]{0.06\linewidth}
\begin{sideways}
\parbox{2.9cm}{\centering Adjoint-Magnitude \\ Flagging Grids}
\end{sideways}
\end{minipage}
\begin{minipage}[c]{0.23\linewidth}
\includegraphics[width=\textwidth]{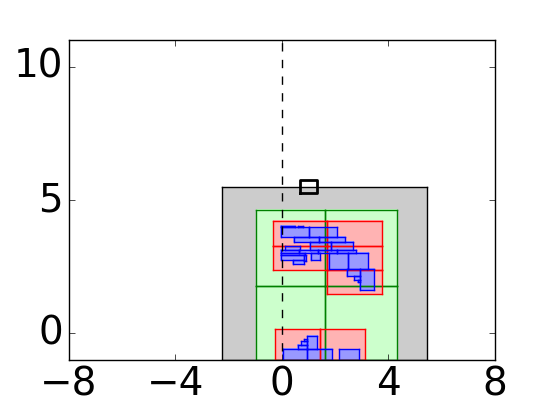}
\end{minipage}
\begin{minipage}[c]{0.23\linewidth}
\includegraphics[width=\textwidth]{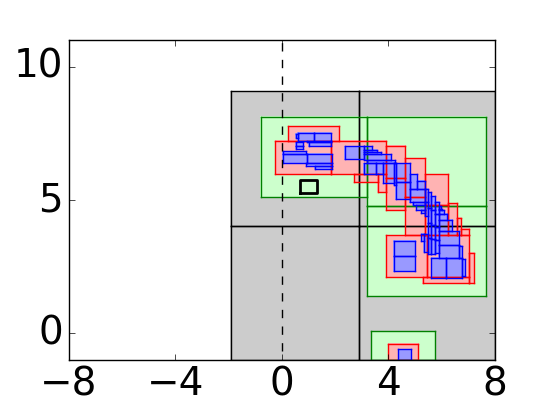}
\end{minipage}
\begin{minipage}[c]{0.23\linewidth}
\includegraphics[width=\textwidth]{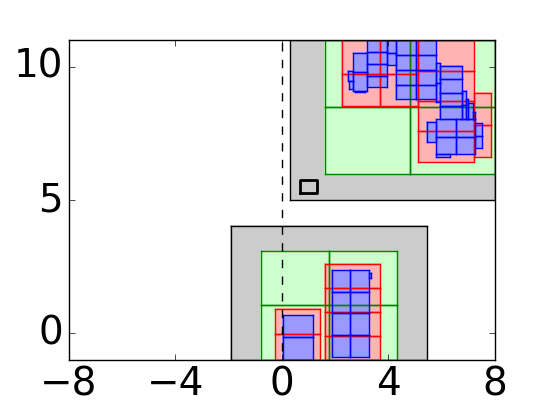}
\end{minipage}
\begin{minipage}[c]{0.23\linewidth}
\includegraphics[width=\textwidth]{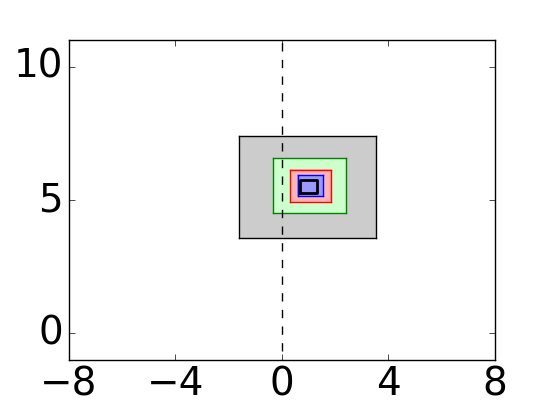}
\end{minipage}\\
\begin{minipage}[c]{0.06\linewidth}
\begin{sideways}
\parbox{2.9cm}{\centering Adjoint-Error \\ Flagging Grids}
\end{sideways}
\end{minipage}
\begin{minipage}[c]{0.23\linewidth}
\includegraphics[width=\textwidth]{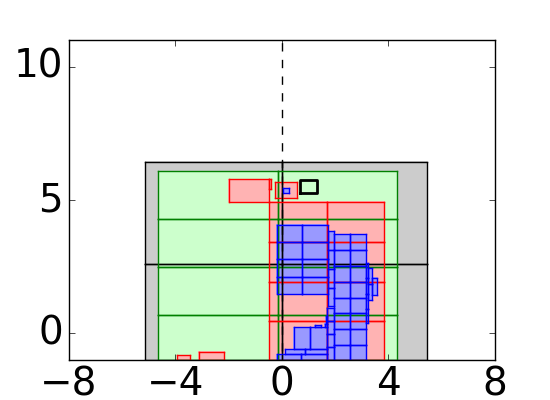}
\end{minipage}
\begin{minipage}[c]{0.23\linewidth}
\includegraphics[width=\textwidth]{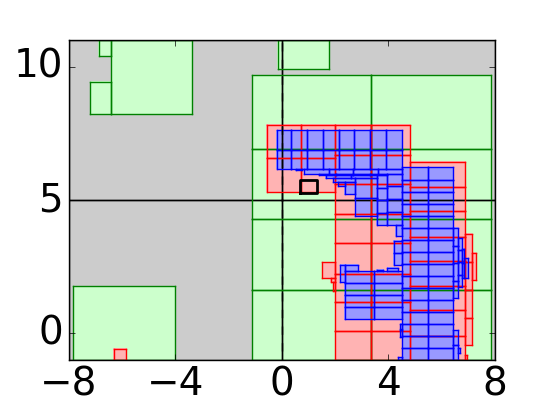}
\end{minipage}
\begin{minipage}[c]{0.23\linewidth}
\includegraphics[width=\textwidth]{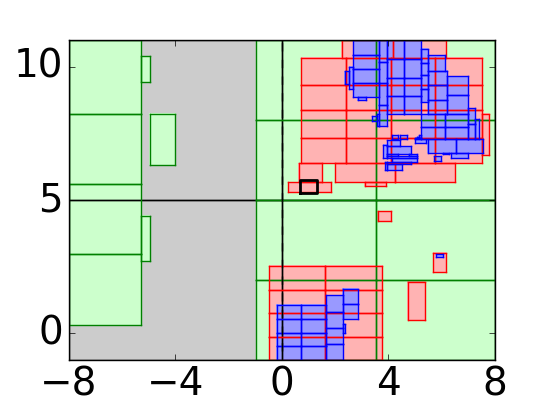}
\end{minipage}
\begin{minipage}[c]{0.23\linewidth}
\includegraphics[width=\textwidth]{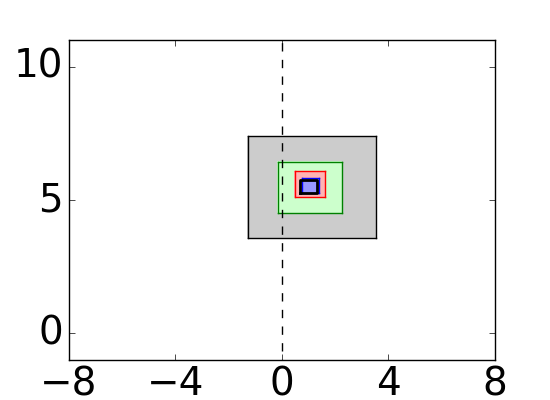}
\end{minipage}
 \caption{%
  Plots for the two dimensional forward problem 
  for example 3. The top row shows 
  the fine grid solution, and the bottom four rows each show the 
  grids for the different refinement levels being used for this 
  problem when it is solved using one of the four flagging methods being 
  presented. 
  The color scale for the fine grid solution figures goes from blue to red, and
  ranges between $-0.3$ and $0.3$. 
  For all of the 
refinement level plots the coarsest refinement level is shown in white, 
and the second, third, fourth and fifth levels of refinement are shown 
in grey, green, red, and blue, respectively. }
   \label{fig:2d_plot_ex3}
\end{figure}

\Cref{fig:2d_plot_ex3} compares the refinement levels 
used by each of the four different flagging methods on the 
forward problem for this example. 
Difference-flagging was used with a tolerance of $10^{-12}$ to 
compute a fine grid solution to our forward problem.
The top row of the figure shows this fine grid solution, 
for the sake of reference. Each of the following rows shows 
the refinement patches (colored by level) 
being used in the simulation, from 
which we can quickly note that the adjoint-flagging methods have smaller
regions of refinement than the other two flagging methods. 
The location of interest is outlined with a black box on all of the plots. 
For each of these flagging methods five levels of refinement are allowed
for the forward problem,
where the coarsest level is a $50 \times 50$ grid 
and a refinement ratio of 2  in both $x$ and $y$ 
is used from each grid to the next. So, the finest level in the forward problem 
corresponds to a fine grid of $1600 \times 1600$ cells. For all of the 
refinement level plots the coarsest refinement level is shown in white, 
and the second, third, fourth and fifth levels of refinement are shown 
in grey, green, red, and blue, respectively. 
Where these finer levels of refinement are placed, of course, varies 
based on the flagging method being used.
The adjoint problem was solved on a relatively 
coarse grid of $200 \times 200$ cells, and 
no adaptive mesh refinement. 
(As a side note, this problem was also computed with the 
adjoint problem solved on a coarse grid of $50 \times 50$ cells 
with nearly the same final results. The $200 \times 200$ 
cells adjoint was used for the figures shown here, since it provided 
clearer plots for \cref{fig:2d_adjoint_ex3}).
The tolerance used for the 
difference-flagging plots shown was $3\times 10^{-2}$, 
for the adjoint-magnitude flagging plots the tolerance used was $3 \times 10^{-4}$, 
for the error-flagging plots the tolerance used was $6\times 10^{-5}$, and 
for the adjoint-error flagging plots the tolerance used was $3\times 10^{-3}$. 
These tolerances were chosen because they all resulted in an error 
of about $5\times 10^{-4}$ in our functional $J$. Therefore, it is 
appropriate to compare the placement of the refined patches
in the plots shown for all four flagging methods. 

Note that, as in the one-dimensional examples, the new adjoint-based approaches
allow refinement to occur only around the waves that eventually coalesce on the
location of interest, with a region very close to the rectangle defined by
the functional $J$ being refined at the final time $t_f=21$.  By contrast, the
original difference-flagging and error-flagging approaches lead to refinement of
many waves and portions of the domain that have no benefit in terms of
determining the functional $J$ accurately.

\subsection{Example 4: Capturing many intersecting waves}
How much of the domain must be refined with the adjoint approach
depends on the specification of $J$.   In this final example we simply move the
location of interest to a place where more waves are converging at $t_f$ in order
to illustrate that the adjoint approach will adapt to this situation.
Suppose that we are now interested in the accurate estimation of the 
pressure in the area defined by a rectangle centered about 
$(x,y) = (3.5,0.5)$, and define
\begin{align*}
J = \int_{3.18}^{3.82}\int_{0.26}^{0.74}p(x,y,t_f)dy\,dx.
\end{align*}

For this example, the definition of the adjoint problem is the same as 
in example 3 except for the ``initial data'', for which we now take
$\hat{q}(x,y,t_f) = \varphi (x,y)$ 
where $\varphi$ is as in \cref{eq:phi_2d_ex3} with $I(x,y)$ now defined by
\begin{align}
I (x,y) = \left\{
     \begin{array}{ll}
       1 & \hspace{0.3in}\textnormal{if } 3.18 \leq x \leq 3.82
        \textnormal{ and } 0.26 \leq y \leq 0.74,\\
       0 &\hspace{0.3in} \textnormal{otherwise.}\hspace{1.61in}
     \end{array}
   \right. \label{eq:delta_2d}
\end{align}
As time progresses backwards, waves radiate outward and reflect off the
walls as well as transmitting and reflecting off of the interface at $x = 0$. 

\begin{figure}[h!]
\begin{minipage}[c]{0.06\linewidth}
\hspace{0.1cm}
\end{minipage}
\begin{minipage}[c]{0.23\linewidth}
\centering
t = 2.5 seconds
\end{minipage}
\begin{minipage}[c]{0.23\linewidth}
\centering
t = 6 seconds
\end{minipage}
\begin{minipage}[c]{0.23\linewidth}
\centering
t = 15 seconds
\end{minipage}
\begin{minipage}[c]{0.23\linewidth}
\centering
t = 21 seconds
\end{minipage}\\
\begin{minipage}[c]{0.06\linewidth}
\begin{sideways}
\parbox{2.9cm}{\centering Fine Grid \\ Solution}
\end{sideways}
\end{minipage}
\begin{minipage}[c]{0.23\linewidth}
\includegraphics[width=\textwidth]{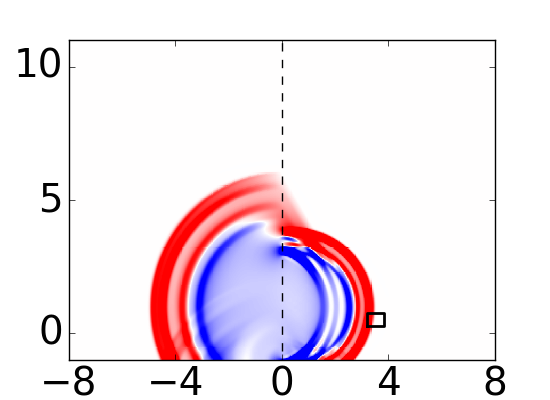}
\end{minipage}
\begin{minipage}[c]{0.23\linewidth}
\includegraphics[width=\textwidth]{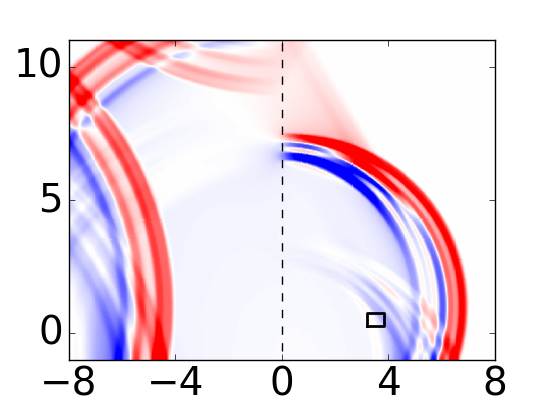}
\end{minipage}
\begin{minipage}[c]{0.23\linewidth}
\includegraphics[width=\textwidth]{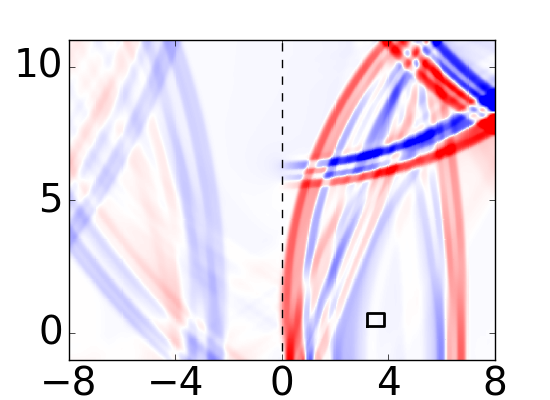}
\end{minipage}
\begin{minipage}[c]{0.23\linewidth}
\includegraphics[width=\textwidth]{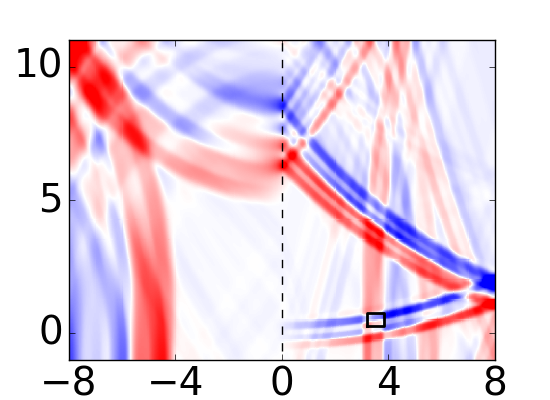}
\end{minipage}\\
\begin{minipage}[c]{0.06\linewidth}
\begin{sideways}
\parbox{2.9cm}{\centering Adjoint-Magnitude \\ Flagging Grids}
\end{sideways}
\end{minipage}
\begin{minipage}[c]{0.23\linewidth}
\includegraphics[width=\textwidth]{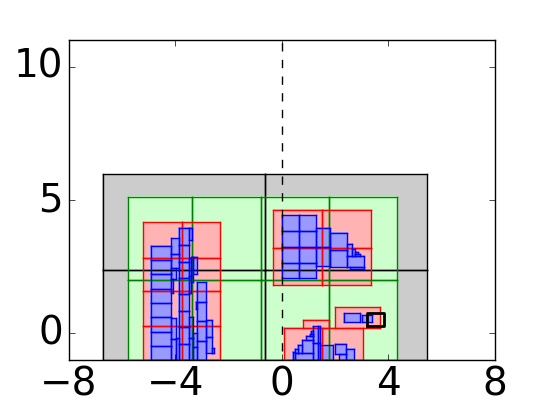}
\end{minipage}
\begin{minipage}[c]{0.23\linewidth}
\includegraphics[width=\textwidth]{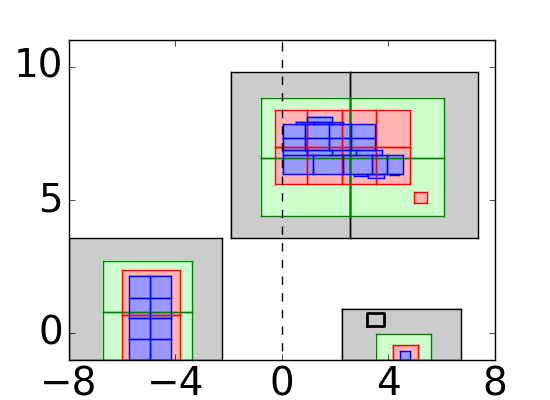}
\end{minipage}
\begin{minipage}[c]{0.23\linewidth}
\includegraphics[width=\textwidth]{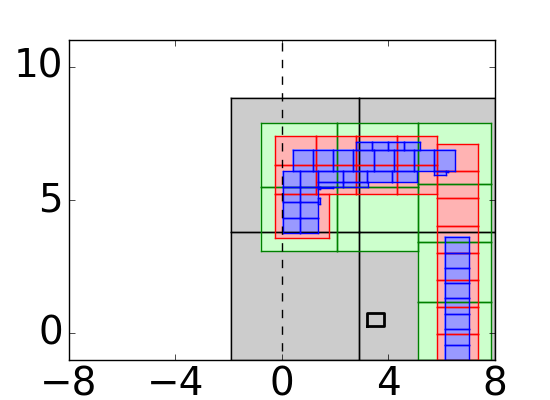}
\end{minipage}
\begin{minipage}[c]{0.23\linewidth}
\includegraphics[width=\textwidth]{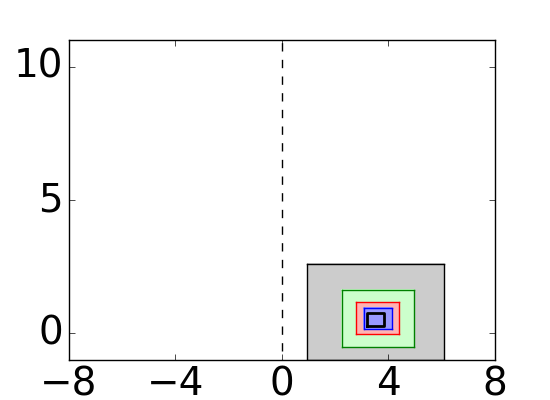}
\end{minipage}\\
\begin{minipage}[c]{0.06\linewidth}
\begin{sideways}
\parbox{2.9cm}{\centering Adjoint-Error \\ Flagging Grids}
\end{sideways}
\end{minipage}
\begin{minipage}[c]{0.23\linewidth}
\includegraphics[width=\textwidth]{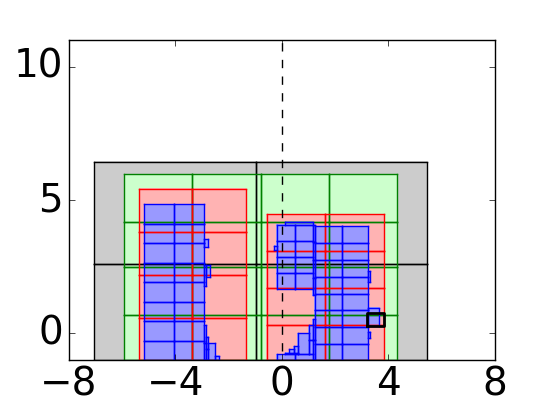}
\end{minipage}
\begin{minipage}[c]{0.23\linewidth}
\includegraphics[width=\textwidth]{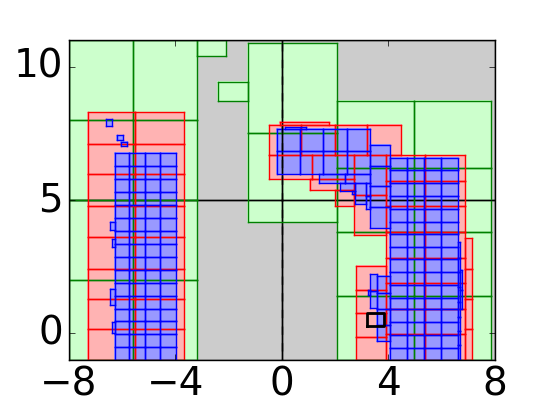}
\end{minipage}
\begin{minipage}[c]{0.23\linewidth}
\includegraphics[width=\textwidth]{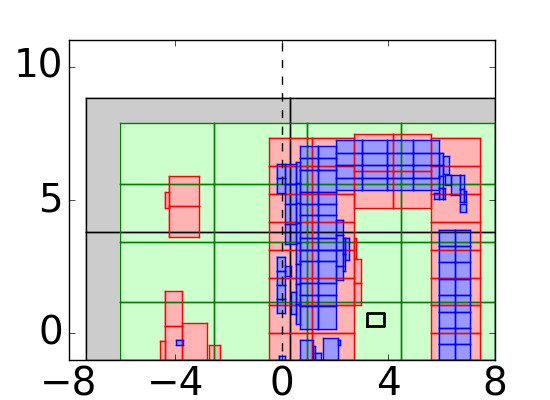}
\end{minipage}
\begin{minipage}[c]{0.23\linewidth}
\includegraphics[width=\textwidth]{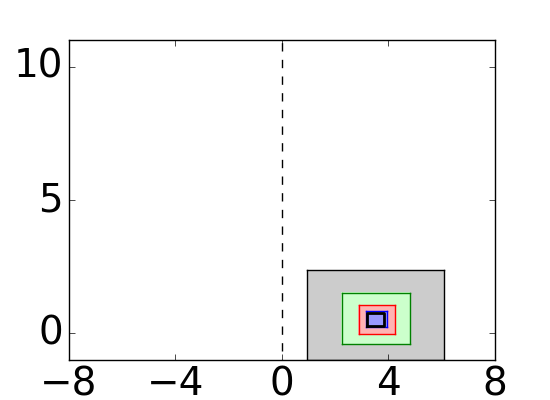}
\end{minipage}
 \caption{%
  Plots for the two dimensional forward problem 
  for example 4. The top row shows 
  the fine grid solution, and the bottom two rows each show the 
  grids for the different refinement levels being used for this 
  problem when it is solved using one of the two adjoint-flagging methods. 
  The color scale for the fine grid solution figures goes from blue to red, and
  ranges between $-0.3$ and $0.3$. 
  For all of the 
refinement level plots the coarsest refinement level is shown in white, 
and the second, third, fourth and fifth levels of refinement are shown 
in grey, green, red, and blue, respectively. }
   \label{fig:2d_plot}
\end{figure}

\Cref{fig:2d_plot} compares the refinement levels 
used by each of the two adjoint-flagging methods on the 
forward problem for this example. 
The top row of the figure shows the fine grid solution once again, 
for the sake of reference and to show the new location of interest. 
The new location of interest is outlined with a black box on all of the plots. 
As with \cref{fig:2d_plot_ex3}, 
each of the following rows shows the refinement patches (colored by level) 
being used in the simulation. 
For each of these flagging methods five levels of refinement are allowed
for the forward problem,
where the coarsest level is a $50 \times 50$ grid 
and a refinement ratio of 2  in both $x$ and $y$ 
is used from each grid to the next. 
So, as with example 3, the finest level in the forward problem 
corresponds to a fine grid of $1600 \times 1600$ cells. For all of the 
refinement level plots the coarsest refinement level is shown in white, 
and the second, third, fourth and fifth levels of refinement are shown 
in grey, green, red, and blue, respectively. 
Where these finer levels of refinement are placed, of course, varies 
based on the flagging method being used.
The adjoint problem was again solved on a relatively 
coarse grid of $200 \times 200$ cells without adaptive mesh refinement. 

The tolerance used for these plots are the same 
as the ones shown for example 3. 
Note that for the same tolerance 
the fine grid, difference flagging, and error flagging 
solutions are the same as in example 3, except for the 
placement of the gauge. Since none of these 
computations depended on the adjoint method, changing the 
area of interest (which changes the initial conditions for 
the adjoint problem) does not affect the computation. 
Since they would be the same as 
\cref{fig:2d_plot_ex3}, the figures for difference-flagging  and 
error-flagging grids are not repeated in 
\cref{fig:2d_plot}.

\subsection{Computational Performance}

\begin{figure}[h!]
\begin{minipage}[b]{0.45\linewidth}
\includegraphics[width=\textwidth]{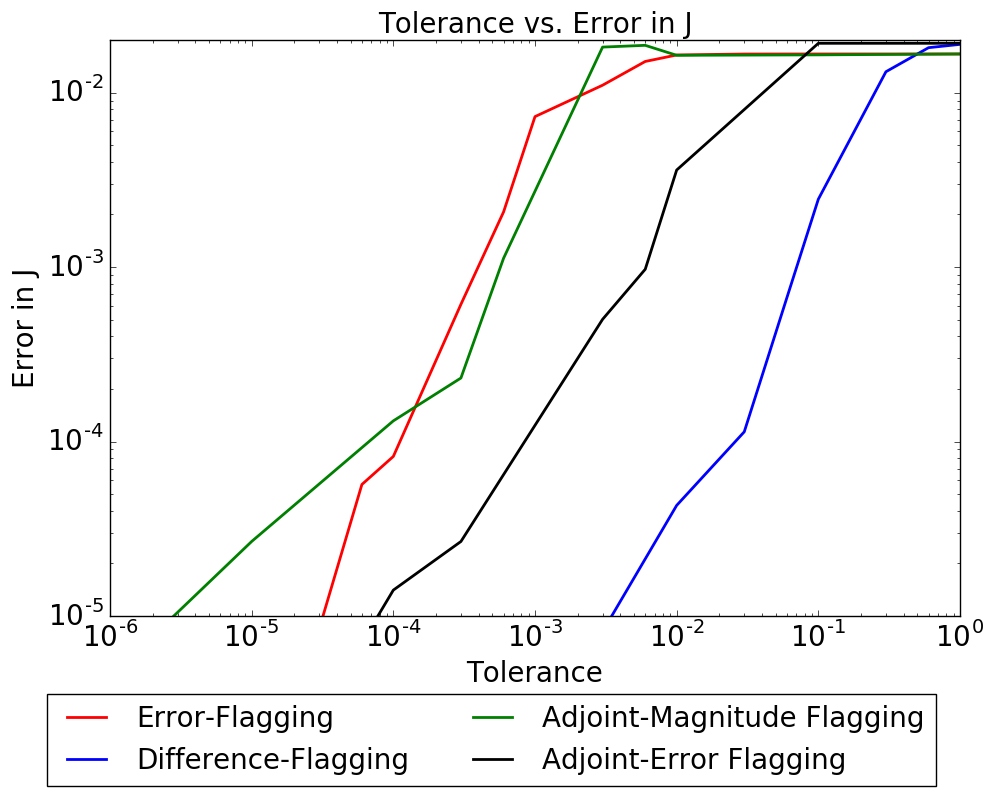}
\end{minipage}
 \caption{%
 Accuracy for the different flagging methods on 
 the two-dimensional acoustics example 3. 
 Shown is the error in the functional $J$ for each tolerance value 
 for the various flagging methods.
 }
\label{fig:2d_tolvserror}
\end{figure}

The advantage of using the adjoint method 
can be seen in the amount of work that is required for 
each of the flagging methods we are considering, and 
the level of accuracy that results from that work. 
\Cref{fig:2d_tolvserror} shows the error in our functional of 
interest for the two-dimensional acoustics example 3 
as the tolerance is varied for the four flagging 
methods we are considering. 
Our fine grid solution was used to calculate a fine grid value of our functional 
of interest, $J_\text{fine}$. 
The error 
between the calculated value of $J$, from the various tolerances used 
for each of the flagging methods, 
and $J_\text{fine}$ is shown. 

Of particular interest, note that the magnitude of the error when using adjoint-error 
flagging is once again consistently less than the magnitude of the tolerance being used, 
as we expected based on \cref{sec:adjErrorFlag}. This 
means that, as with example 1, adjoint-error flagging allows the user to enforce 
a certain level of accuracy on the final functional of interest by 
selecting a tolerance of the desired order. 
Recall that the tolerance set by the user is being used 
to evaluate whether or not some given quantity is above that set 
tolerance. However, the quantity that is being evaluated is different 
for each of the flagging methods we are considering. Therefore, 
other than noting the general trends in this figure, not 
much benefit comes from comparing each of the lines 
corresponding to a flagging method to the others in this plot. 

\begin{figure}[h!]
\begin{minipage}[b]{0.45\linewidth}
\includegraphics[width=\textwidth]{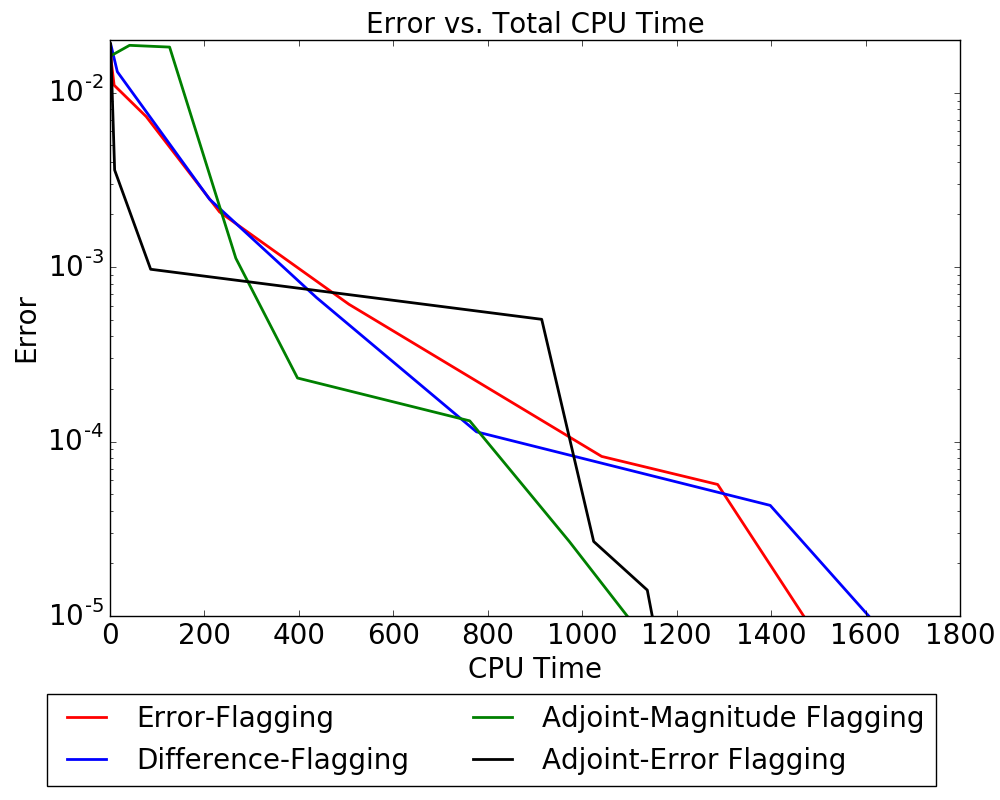}
\end{minipage}
\vspace{0.2cm}
\begin{minipage}[b]{0.45\linewidth}
\includegraphics[width=\textwidth]{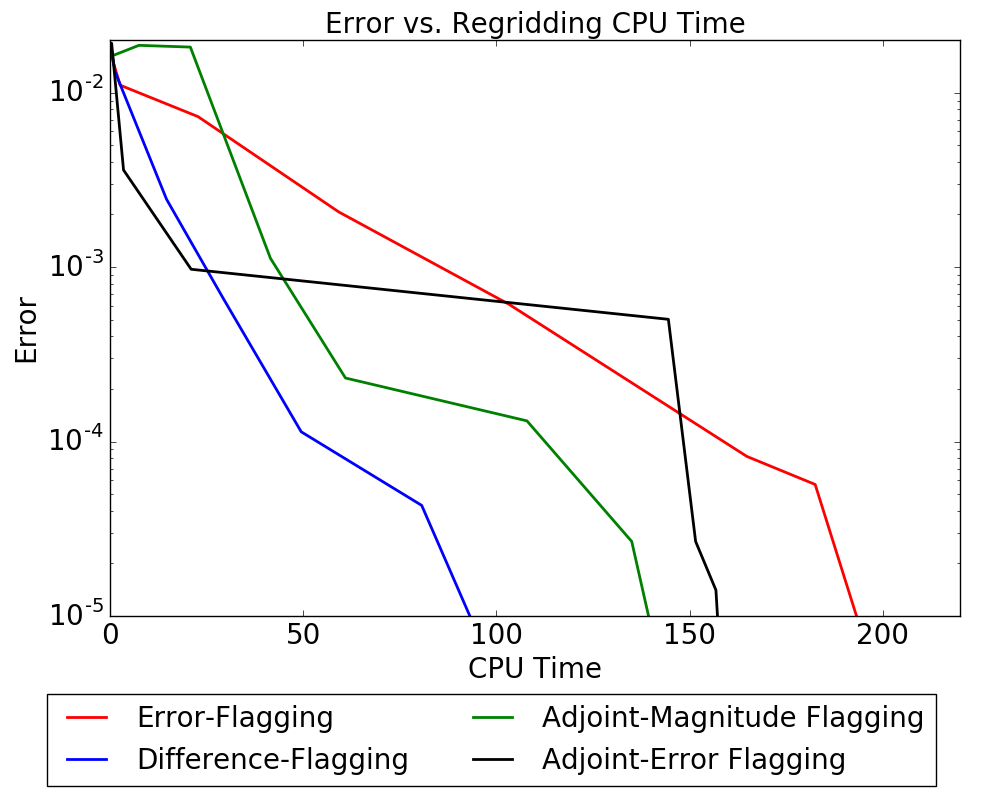}
\end{minipage}
 \caption{%
Performance measures for the different flagging methods on 
 the two-dimensional acoustics example 3. 
 On the left: the total CPU time (in seconds) required vs. the 
 accuracy achieved.
 On the right: the regridding CPU time (in seconds) required vs. the 
 accuracy achieved.  CPU times were found by averaging the CPU time 
 over ten runs.}
\label{fig:2d_accuracy}
\end{figure}

We are better able to compare the methods to one another if we 
consider the amount of computational time that is required to 
achieve a certain level of accuracy with each flagging method. 
As with the one-dimensional examples, 
both of the two-dimensional examples was run on a quad-core laptop, 
and the OpenMP option of AMRClaw was enabled which 
allowed all four cores to be utilized. \Cref{fig:2d_accuracy} 
shows two measures of the amount of CPU time that was 
required for each method vs. the accuracy achieved for 
example 3. On the left, 
we have shown the total amount of CPU time used by the 
computation. This includes stepping the solution forward in time 
by updating cell values, 
regridding at appropriate time intervals, outputting the results, 
and other various overhead requirements. 
For each flagging 
method and tolerance this example was run ten times, 
and the average of the CPU times for those runs was used 
to generate this plot.  
As with the one-dimensional example, while 
the adjoint-magnitude  and adjoint-error flagging  
methods 
have a higher CPU time requirement for lower accuracy, 
these methods quickly showed their strength by maintaining 
a low CPU time requirement while increasing the accuracy of the 
solution. In contrast, the CPU time requirements for 
difference-flagging and error-flagging 
are larger than their adjoint-flagging counterparts when increased 
accuracy is required. 
Note that the adjoint-flagging methods do have the 
additional time requirement of solving the adjoint 
problem, which is not shown in these figures. For this 
example, solving the adjoint required about 37 seconds of CPU 
time, which is once again small compared to the time spent 
on the forward problem. As before, this CPU time was found by taking the average 
time required over ten simulations. 

\begin{figure}[h!]
\begin{minipage}[b]{0.45\linewidth}
\includegraphics[width=\textwidth]{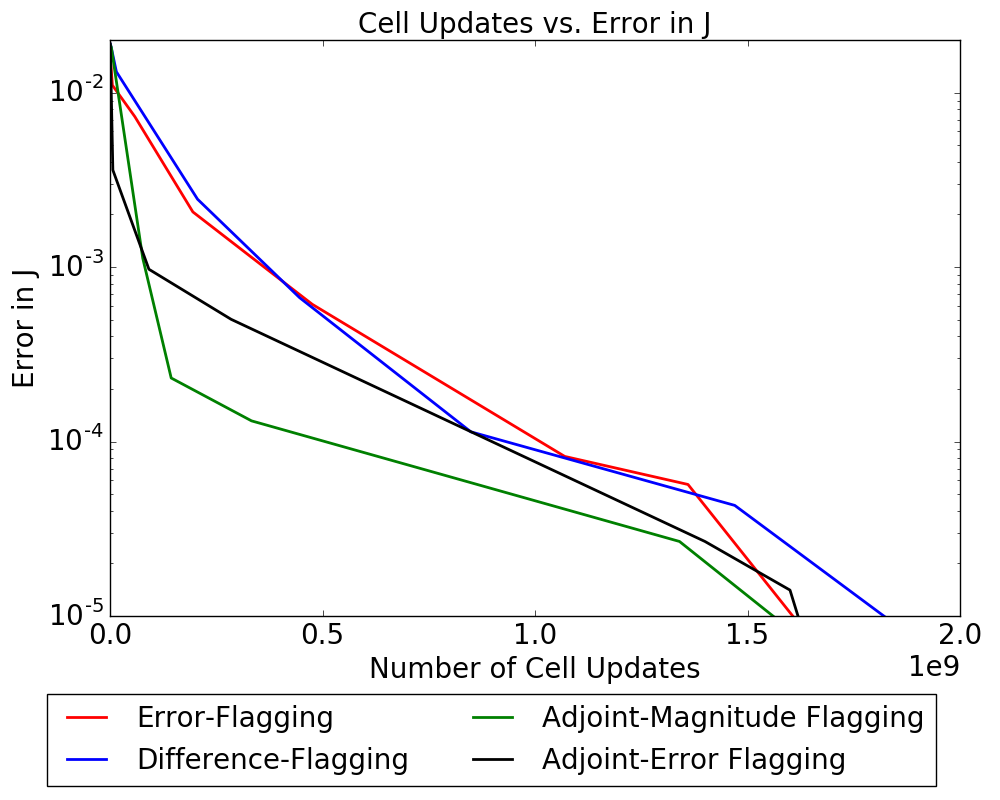}
\end{minipage}
 \caption{%
Number of cell updates calculated vs. the accuracy achieved 
for the four different flagging methods begin considered for 
example 3. Note 
that the number of cell updates axis is multiplied by $1$e$9$.
 }
\label{fig:ex3_cellsvserror}
\end{figure}

The right of \cref{fig:2d_accuracy} shows the amount of CPU time 
spent on the regridding process vs. the accuracy of the functional $J$. 
Note that this is really where the differences between the four methods 
lie: the adjoint-flagging methods goal is to reduce the number of cells 
that are flagged while maintaining the accuracy of the solution. While 
we have seen from the left side of \cref{fig:2d_accuracy} that 
this has certainly been successful 
for adjoint-flagging in reducing the overall time required 
for computing the solution, we might expect that the time spent in the 
actual regridding process might be longer for the adjoint-flagging methods. 
 However the amount of time spent in the 
regridding process is not significantly greater for the adjoint-flagging methods. 

As with the one-dimensional example, here we see that the adjoint-flagging 
methods have less memory requirements than the alternative flagging methods. 
\Cref{fig:ex3_cellsvserror} shows 
the number of cell updates required vs. the accuracy achieved for each of the 
flagging methods being considered. While the results 
are not as drastic as what we saw for example 1, note that the number of cell updates 
for example 3 is consistently less for the adjoint-flagging methods than for their non-adjoint 
counterparts, since there are fewer fine grids 
throughout the domain. This reduces the memory requirements for the computation. 
This example is also relatively small, and memory usage when utilizing the non-adjoint 
flagging method did not become an issue. However, in larger simulations (for instance, 
ocean-wide problems when considering tsunami modeling) memory constraints can 
become a significant consideration. It can also become a significant consideration 
in three dimensional problems. 

\section{Conclusions and Future Work}

In this paper we have presented a method for using the adjoint equation 
to guide adaptive mesh refinement when solving hyperbolic systems 
of equations, as well as describing the implementation of this method 
into the AMRClaw software package. The adjoint method 
was used to identify the waves which should be refined to ensure 
an accurate estimation of the solution at the final time, provided 
we were interested in some target region (for example, the location 
of a pressure gauge). 

Two different methods for identifying cells 
that should be refined were presented: adjoint-magnitude flagging 
and adjoint-error flagging. These two methods were compared 
to the two main flagging methods currently available in AMRCLAW: 
difference-flagging and error-flagging. 
Accuracy and timing results for 
various examples were presented, and for 
each example the results from using 
all four of these flagging methods were compared. 

It was shown 
that using adjoint-flagging provides confidence that the 
appropriate waves have been refined to accurately capture the 
solution at the target location, without unnecessary refinemement in other
regions. For the examples considered, both adjoint-magnitude 
flagging and adjoint-error flagging had the advantage of requiring less 
computational time than their non-adjoint counterparts to 
achieve the same level of accuracy in the final solution. Also, 
the adjoint-flagging methods required fewer cell updates to 
achieve the same level of accuracy, which signifies a smaller 
memory requirement. 

The code for all the examples presented in this work is available 
online at \cite{adjointCode}. This code can be easily modified to 
solve other linear variable coefficient acoustics problems. This 
repository also contains other examples illustrating how 
adjoint-flagging can be used with AMRClaw and GeoClaw.
These adjoint-flagging procedures are now being incorporated into the
Clawpack software and should be available in future releases for general use.

In this work we have focused on examples involving the linear 
acoustics equations. However, the adjoint method is not 
limited to working on only these equations. For examples 
using the adjoint method to guide adaptive mesh refinement 
for the shallow water equations in the context of tsunami 
modeling see \cite{DavisLeVeque:adjoint2016, Borrero2015}. 
Using the adjoint method when solving non-linear equations 
requires that the equations be linearized. In the context of using the
shallow water equations for tsunami propagation across the ocean
this is fairly straight forward, since 
they can be linearized about the ocean at rest. An area of 
future work is studying the use of the adjoint method in the 
context of other nonlinear hyperbolic equations, where the adjoint 
equation would be derived by linearizing about a particular 
forward solution. This would require the development of some 
kind of automated process to shift between solving 
the forward problem, linearizing about that forward problem, 
solving the corresponding adjoint problem, and using that adjoint 
solution in guiding the adaptive mesh refinement for the 
forward problem. 

We have focused 
in this work on the accurate calculation of the functional $J$, where 
the definition for this functional varied based on the target area of
interest.  However, we are typically concerned with the 
accurate estimation of the solution to the forward problem rather than 
the functional $J$. For the examples presented in this work we 
used various different initial conditions for the adjoint problem, but did 
not consider the effect this had on the accuracy of the forward solution. 
Examining the implications of the initial conditions used for the 
adjoint problem, and their effect on the accuracy of the foward 
solution,  is another area for future work. 

\appendix
\section{Riemann solvers used in this work}

In this appendix, we present some details regarding the 
Riemann problems encountered when solving the 
forward and adjoint equations for variable-coefficient linear acoustics. 
A Riemann problem is the hyperbolic equation of interest together with initial 
data that is piecewise constant with a single jump discontinuity, say, 
\begin{align*}
&q_0(x) = \left\{
     \begin{array}{lr}
       q_l & \textnormal{if } x < 0\\
       q_r & \textnormal{if } x > 0
     \end{array}
   \right.
\end{align*}
where $q_l$ and $q_r$ are constant state vectors. 
For variable-coefficient linear problems,
the Riemann problem also assumes that the variable
coefficients have a jump discontinuity at $x=0$, e.g., from $\rho_l$ to
$\rho_r$.  
  
For the acoustics equations $q_t + A(x)q_x=0$ with $A(x)$ given by
\cref{eq:Aacoustics}, it is easy to compute that the eigenvalues are 
\begin{align*}
\lambda^1 = -c(x) \hspace{0.3in} \textnormal{and} \hspace{0.3in}\lambda^2 = c(x)
\end{align*}
where $c(x) = \sqrt{K(x)/\rho(x)}$ is the speed of sound in material. The 
eigenvectors for this matrix $A(x)$ are 
\begin{align*}
r^1 =  \left[\begin{matrix}
-Z(x) \\ 1
\end{matrix}\right]
\hspace{0.3in} \textnormal{and} \hspace{0.3in}
r^2 = \left[\begin{matrix}
Z(x) \\ 1
\end{matrix}\right],
\end{align*}
where $Z(x) = \rho(x)c(x)$ is the impedance.
Note that any scalar multiple of these vectors would still be an 
eigenvector of $A(x)$. This particular normalization was chosen 
for consistency with \cite{Leveque1}, where many more details are given about
Riemann solutions and the manner in which these are used as the building block
for the high-resolution finite volume methods implemented in Clawpack.

\begin{figure}[h!]
\begin{minipage}[b]{0.5\linewidth}
\includegraphics[width=\textwidth]{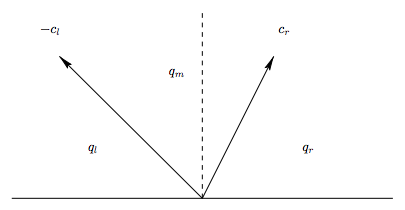}
\end{minipage}
 \caption{%
Structure of the solution to the Riemann problem for the forward problem, 
in the $x-t$ plane. The dashed 
line depicts the interface between two different cells. 
Between the waves is a single state $q_m$. 
For the adjoint problem with initial states $\tilde q_l$ and $\tilde q_r$
the wave speeds are the same but there will be two intermediate
states $\tilde q_{ml}$ and $\tilde q_{mr}$ to either side of the dashed line.
Figure from \cite{Leveque1}.}
\label{fig:q_waves}
\end{figure}

The solution to the Riemann problem for acoustics consists of two waves, 
a left-going sound wave with speed $-c_l$ and a right-going sound wave with
speed $c_r$.  The left-going wave moves into a homogeneous material with
impedance $Z_l$ and hence the jump in $q$ across this wave is a multiple of the
eigenvector $r^1_l = [-Z_l, 1]^T$, while the jump in $q$ across the right-going
wave must be a multiple of $r^1_r = [Z_r, 1]^T$.  Between these two waves the
state $q_m$ must be constant across the jump in material properties at $x=0$, or
else the resulting stationary jump in $q$ would lead to a delta function
singularity in $q_t$ according to the PDE.  This structure is illustrated in
\cref{fig:q_waves}.

Therefore, we wish to find $q_m$ such that
 \begin{align*}
 q_m - q_l = \alpha^1\left[\begin{matrix}
-Z_l \\ 1
\end{matrix}\right]  \equiv \mathcal{W}^1_{rl} 
\hspace{0.3in} \textnormal{and} \hspace{0.3in}
q_r - q_m = \alpha^2 \left[\begin{matrix}
Z_r \\ 1
\end{matrix}\right] \equiv \mathcal{W}^2_{rl} 
 \end{align*}
for some scalar coefficients $\alpha^1$ and $\alpha^2$. 
If we add these two equations together we get 
\begin{align*}
q_r - q_l = \alpha^1\left[\begin{matrix}
-Z_l \\ 1
\end{matrix}\right] + \alpha^2 \left[\begin{matrix}
Z_r \\ 1
\end{matrix}\right].
\end{align*}
This gives us a linear system of two equations to solve for $\alpha^1$ 
and $\alpha^2$. 
The solution of this gives us the decomposition of the jump in $q$ as the sum of
the two acoustic waves,
\begin{equation}\label{eq:wavesplitting}
  q_r - q_l = \sum\limits_{p = 1}^2 \alpha^p r^p \equiv \sum\limits_{p =1}^2\mathcal{W}^p_{rl}
\end{equation}

When solving the Riemann problem for the modified adjoint equation
$\tilde q_t - (A(x)^T\tilde q)_x = 0$
that is solved forward in time (see  \cref{eq:1d_adjoint_modified}),
the eigenvalues of $-A(x)^T$ are the same as the eigenvalues of $A(x)$, but the
eigenvectors are now
\begin{align*}
\tilde{r}^1 =  \left[\begin{matrix}
1 \\ Z(x)
\end{matrix}\right]
\hspace{0.3in} \textnormal{and} \hspace{0.3in}
\tilde{r}^2 = \left[\begin{matrix}
1 \\ -Z(x)
\end{matrix}\right].
\end{align*}
Since the adjoint equation is in conservation form, it is the flux 
$\tilde f(\tilde q,x) \equiv -A(x)^T\tilde q$ that must be continuous across $x=0$ in order
to avoid a singularity, and so $\tilde q$ will generally have have a jump from
$\tilde q_{ml}$ to $\tilde q_{mr}$ across $x=0$, 
while $\tilde f(\tilde q_{ml}) = \tilde f(\tilde q_{mr}) \equiv \tilde f_m$.

As before, the jump across each wave is given by 
an eigenvector of the coefficient matrix from the appropriate material.
Therefore, we have that 
\begin{align*}
\tilde{f}_m - \tilde{f}_l = \beta^1\left[\begin{matrix}
1 \\ Z_l
\end{matrix}\right] 
\hspace{0.3in} \textnormal{and} \hspace{0.3in}
\tilde{f}_r - \tilde{f}_m = \beta^2 \left[\begin{matrix}
1 \\ -Z_r
\end{matrix}\right] 
 \end{align*}
for some scalar coefficients $\beta^1$ and $\beta^2$. 
If we add these two equations together we get 
\begin{align*}
\tilde{f}_r - \tilde{f}_l = \beta^1\left[\begin{matrix}
1 \\ Z_l
\end{matrix}\right] + \beta^2 \left[\begin{matrix}
1 \\ -Z_r
\end{matrix}\right].
\end{align*}
This gives us a linear system of equations to solve for $\beta^1$ 
and $\beta^2$.  These so-called f-waves 
can be used
directly in the f-wave version of the wave-propagation algorithms used in
Clawpack, as discussed further in \cite{BaleLevMitRoss02}, or these can be
converted into jumps in $\tilde q$ as follows.
Since the left-going and right-going waves propagate through constant materials,
we have (by the Rankine-Hugoniot jump conditions) that
\begin{align*}
\tilde q_{ml} - \tilde q_{l} = \frac{\tilde f_m - \tilde f(\tilde
q_l)}{\lambda^1_l}
= -\frac{\beta^1}{c_l} \left[\begin{matrix} 1 \\ Z_l \end{matrix}\right],
\qquad
\tilde q_{r} - \tilde q_{ml} = \frac{\tilde f(\tilde q_r)) - \tilde
f_m}{\lambda^2_r}
= \frac{\beta^2}{c_r} \left[\begin{matrix} 1 \\ Z_r \end{matrix}\right],
\end{align*}
which allows determining the jumps in $\tilde q$ across each propagating wave.

More complete details can be found in the references cited above, and
in the Clawpack code implementing the examples shown in this paper,
all of which can be found in \cite{adjointCode}.

\bibliographystyle{ACM-Reference-Format}
\bibliography{toms_bib}